\input amstex
\documentstyle{amsppt}

\def\bbR{{\Bbb R}}
\def\bbZ{{\Bbb Z}}

\def\bbC{{\Bbb C}}

\def\bbP{{\Bbb P}}

\def\bbF{{\Bbb F}}

\def\calM{{\Cal M}}

\def\calO{{\Cal O}}

\def\calX{{\Cal X}}
\def\calV{{\Cal V}}
\def\calI{{\Cal I}}
\def\calC{{\Cal C}}

\hsize 6in
\topmatter
\title 
INVARIANT STABLE BUNDLES OVER MODULAR
CURVES
${\bold X}({\bold p})$
\endtitle

\author 
Igor V. Dolgachev
\endauthor
\thanks 
Research supported in part by NSF grant  DMS
9623041
\endthanks
\endtopmatter

\document

\vglue 1in
\noindent
{\bf Introduction.}

Let $X$ be a smooth projective algebraic curve of genus $g > 1$
and $G$ be the group of its automorphisms. The problem is to describe vector bundles
on $X$ which are invariant with respect to the action of $G$ on $X$. In this paper we
address this problem in the case when the curve $X$ is the modular curve $X(p)$
obtained as a compactification of the quotient of the upper-half plane $H = \{z\in
\bbC:\Im z > 0\}$ by the action of the principal congruence subgroup $\Gamma(p) = \{A =
\pmatrix a&b\\
 c&d\endpmatrix\in SL(2,\bbZ):A\equiv I \quad\hbox{mod $p$}\}$. We shall
assume that $p$ is a prime number $ > 5$ although some of our results are true for any
$p$ not divisible by $2$ and $3$.  Here the group
$G$ is isomorphic to the group
$PSL(2,\bbF_p)$. Also we restrict ourselves with stable bundles. In other words, we are
trying to describe the set of fixed points for the natural action of $G$ on the moduli
space
 of rank $r$ stable vector bundles on $X(p)$. The case of rank 1
bundles is rather easy and the answer can be found in {\bf [AR]}. The group of
$G$-invariant line  bundles on $X(p)$ is generated by a line bundle $\lambda$ of
degree ${p^2-1\over 24}$ which is a $(2p-12)$-th root of the canonical bundle. For
the future use we generalize this result to any Riemann surface $X$ with a finite
group $G$ of its automorphims such that $X/G\cong \bbP^1$. This result must be known 
to experts but I could not find a reference. 
When the determinant of the bundle is trivial, we are able to relate our problem to the
problem of classifying  unitary representations of the fundamental group of the
Brieskorn  sphere
$\Sigma(2,3,p)$, that is, the link of the singularity 
$x^2+y^3+z^p =0$. Applying some known results from differential
topology we prove that there exist exactly
$2n$ rank 2 $G$-invariant stable bundles with trivial determinant and $3n^2\pm n$ rank 3 (if $p \ne 7)$
$G$-invariant stable bundles with trivial determinant on
$X(p)$, where $p = 6n \pm 1$. 
Note that the determinant of a stable $G$-invariant rank 2 bundle is an even multiple of $\lambda$.
So, after twisting by a line bundle, we obtain a $G$-invariant bundle with trivial determinant.  
 
Even in the case of rank $2$ and $3$ our results are still unsatisfactory since we
were able to give a geometric construction of these bundles only in the case $p = 7$.
Some of the bundles we discuss are intrinsically related to the beautiful geometry of
modular curves which goes back to Felix Klein.
	
	I would like to thank the organizers of the conference for giving me the opportunity to
revisit Korea. This paper owes much to the work of Allan Adler and correspondence with him. The book {\bf
[AR]} was a great inspiration for writing this paper.
Finally I would like to thank Hans Boden and Nikolai Saveliev for coaching me in
the theory of Casson invariant of 3-dimensional manifolds.

\bigskip
\flushpar{\bf 1. $G$-invariant and $G$-linearized stable vector bundles.} Let $X$ be a
compact Riemann surface of genus
$g$. For each $r > 0$ there is the moduli space
$\calM_X(r)$ of semi-stable rank $r$ vector bundles
over $X$. Assume that
a finite group $G$ acts holomorphically on $X$
(not necessary faithfully). By functoriality $G$
acts holomorphically on $\calM_X(r)$ and we denote by
$\calM_X(r)^G$ the subvariety of fixed points of
this action. If $[E]\in \calM_X(r)$ is the isomorphism class of
a stable bundle, then $[E]\in \calM_X(r)^G$ if and only if $E$
is $G$-invariant, i.e. for any
$g\in G$, there is an isomorphism of vector bundles 
$$\phi_g:g^*(E)\to E.$$
If $[E]\in \calM_X(r)$ is the point representing the
equivalence class of a semi-stable but not stable
bundle $E$, then it is known that $E$ is equivalent to a
decomposable bundle $E' = E_1\oplus\ldots\oplus E_k$, where all
$E_i$ are stable of the same slope $\mu(E_i) = {deg(E_i)\over rk E_i}$ as $E$ (see {\bf [Se]}). Then
$[E]\in \calM_X(r)^G$ if and only if $E'$ is $G$-invariant.
In the following we will always assume that $E$ is either
stable or is decomposable as above. Assume that the collection
$\phi =
\{\phi_g\}_{g\in G}$ can be chosen in such a way
that for any
$g,g'\in G$
$$\phi_{g\circ g'} = \phi_{g'}\circ g'{}^*(\phi_g).$$
Then we say that $E$ admits a $G$-linearization,
and the pair $(E,\phi)$ is called a
{\it $G$-linearized vector bundle}. Of course, in
down-to-earth terms this means that the action  of
$G$ on $X$ lifts to an action on the total space
of $E$ which is linear on each fibre and a $G$-linearization is such a lift.
One naturally defines the notion of a morphism of
$G$-linearized vector bundles, and, in
particular, one defines the set $\calM_X(G;r)$ of
isomorphism classes of $G$-linearized semi-stable
rank $r$ vector bundles over $G$. There is a
natural forgetting map
$$e: \calM_X(G;r)\to \calM_X(r)^G.$$

\medskip
\plainproclaim 
Proposition 1.1.  Let $E$ be a stable $G$-invariant rank $r$ bundle on $X$. 
One can assign to $E$ an
element
$$c(E)\in H^2(G,\bbC^*)$$
such that $E$ belongs to the image of
the map $e$ if and only if $c(E) = 1$. Here the
cohomology group is taken with respect to the
trivial action of $G$ on the group $\bbC^*$.
\par

{\sl Proof.} This is of course rather standard.
Let $\phi_g:g^*(E)\to E, g\in G,$ be some set
of isomorphisms defined by $E$. We have
$$\phi_{gg'} = c_{g,g'}\phi_{g'}\circ
g'{}^*(\phi_{g})$$
for some $c_{g,g'}\in Aut(E)$. Since $E$ is stable, $Aut(E)$ consists only of
homotheties, so that $Aut(E) =
\bbC^*$. It is easy to check that $\{c_{g,g'}\}_{g,g'\in G}$ defines a
2-cocycle of $G$ with coefficients in the
group $\bbC^*$. Its cohomology class $c(E)$ does
not depend on the choice of $\{\phi_g\}_{g\in G}$. It is
trivial if and only if 
$c_{g,g'} = c_{g'}\circ g'{}^*(c_{g})$ for some map
$c:G\to Aut(E), g\to c_g$. Replacing $\phi_g$
with 
$\psi_g = c_g\circ\phi_g$, we get
$$\psi_{gg'} = \psi_{g'}\circ g'{}^*(\psi_g).$$ 
The set $\{\psi_g\}_{g\in G}$ defines a
$G$-linearization on $E$. Clearly for any $E$ in
the image of $e$, we have $c(E) = 1$. This
checks the assertion.

\medskip
\plainproclaim Corollary 1.2. Assume $G$ is a perfect group (i.e. coincides with its
commutator subgroup). Let
$$ 1\to H^2(G,\bbC^*) \to \tilde G \to G\to 1$$
be the universal central extension of the group $G$
defined by the group of Schur multipliers
$H^2(G,\bbC^*)$. Consider the action of $\tilde G$
on $X$ defined by the action of $G$ on $X$. Then
each stable $G$-invariant bundle $E$ admits a
$\tilde G$-linearization.

{\sl Proof.} Use that 
$$H^2(\tilde G,\bbC^*) = 1.$$

\bigskip
Now let us describe the fibres of the map $e$. 
Let $(E,\phi)$ be a $G$-linearized semi-stable
bundle. Consider the direct product $G\times Aut(E)$ which
acts on $X$ via the action of $G$ on $X$ and
the trivial action of $Aut(E)$ on $X$. Obviously
$$E\in \calM_X(r)^{G\times Aut(E)}.$$

\medskip
\plainproclaim Proposition 1.3. Let $E$ be a $G$-linearized vector bundle and let 
$\{\phi_g\}_{g\in G}$ be the family of isomorphisms $\phi_g:g^*(E)\to E$ defining its
linearization.  For each
$(g,\alpha)\in G\times Aut(E)$ set
$$\phi_{(g,\alpha)} = \alpha^{-1}\circ\phi_g:
(g,\alpha)^*(E) = g^*(E)\to E\to E.$$
Then the set of
isomorphisms of vector bundles
$\phi_{(g,\alpha)}$ defines  a
$G\times Aut(E)$-linearization of $E$ if and only
if, for any
$g\in G$ and any $\alpha\in Aut(E)$,
$$\phi_g\circ g^*(\alpha) = \alpha\circ
\phi_g.\eqno (1.1)$$

{\sl Proof.} It is immediately verified that
$$\phi_{(g',\alpha')}\circ g'{}^*(\phi_{(g,\alpha)})
= \alpha'{}^{-1}\circ\phi_{g'}\circ
g'{}^*(\alpha^{-1})\circ g'{}^*(\phi_g).$$
This is equal to 
$$\phi_{(gg',\alpha\circ\alpha')}
=
\alpha'{}^{-1}\circ\alpha^{-1}\circ
\phi_{gg'} =
\alpha'{}^{-1}\circ\alpha^{-1}\circ\phi_{g'}\circ
g'{}^*(\phi_g)$$
if and only if,
for any $g',\alpha$, 
$$\alpha^{-1}\circ\phi_{g'} = \phi_{g'}\circ
g'{}^*(\alpha^{-1}).$$
This is of course equivalent to the assertion of
the proposition.

\bigskip\noindent
{\bf Definition} A $G$-linearization $\phi = \{\phi_g\}_{g\in G}$
on $E\in \calM_X(G;r)$ is called  {\it distinguished}
if $\phi$ satisfies the
condition $(1.1)$ from Proposition 1.3.

\bigskip
\plainproclaim Theorem 1.4. Let $(E,\phi)\in
\calM_X(G;r)$ be a $G$-linearized
bundle with distinguished linearization. Then any
$G$-linearization $\psi$ on $E$ is equal to
$$\psi_g = \lambda(g^{-1})\circ\phi_g,$$
where 
$$\lambda: G\to Aut(E)$$
is a homomorphism of groups.

{\sl Proof.} First we check that for any
homomorphism of groups $\lambda:G\to Aut(E)$ the
collection $\{\psi_g =
\phi_g\circ\lambda(g)\}_{g\in G}$ defines a
$G$-linearization of $E$.  This is straightforward:
$$\psi_{gg'} =  \lambda((gg')^{-1})\phi_{gg'} =
\lambda(g'{}^{-1})\circ\lambda(g^{-1})\circ\phi_{g'}\circ
g'{}^*(\phi_g) = $$
$$\lambda(g'{}^{-1})\circ\phi_{g'}\circ
g'{}^*(\lambda(g^{-1}))\circ g'{}^*(\phi_g) = \psi_{g'}\circ
g'{}^*(\psi_g).$$
So it is checked. Now suppose we have another
$G$-linearization $\psi = \{\psi_g\}_{g\in G}$ on $E$. Then 
$$\phi_g\circ\psi_g^{-1}:E\to g^*(E)\to E$$
is an automorphism
of $E$. Thus $\psi_g = \lambda(g^{-1})\circ \phi_g$ for some
automorphism $\lambda(g^{-1})$ of $E$. Now reversing the
previous computations we check that the map $ g\to
\lambda(g)$ is a homomorphism of groups.

\medskip
\noindent
{\bf Example 1.5}. Every $G$-linearization on a stable
bundle is distinguished and there is a natural
bijection
$$e^{-1}(e((E,\phi))) \to Hom(G,\bbC^*).$$

\medskip
\noindent
{\bf Example 1.6}. Let $E = {\calO}_X^r$ be the trivial bundle.
It is semi-stable but not stable. Consider the {\it trivial
$G$-linearization} on
$E$ by setting for each $(x,v)\in g^*(E)_x = E_{g\cdot x}$
$$\phi_g(x,v) = (g\cdot x,v).$$
Obviously it is distinguished, and 
$$ e^{-1}({\calO}_X^r) = Hom(G,GL(r,\bbC)).$$

\bigskip
\noindent
{\bf 2. Line bundles}. Let us consider the special case of line bundles. Here we can say
much more. First let us denote by $Pic(X)^G$ the groups of
$G$-invariant line bundles and by $Pic(G;X)$ the group of $G$-linearized line bundles.
The latter group has the following simple interpretation in terms of $G$-invariant
divisors: 

\medskip
\plainproclaim Proposition 2.1. The group $Pic(G;X)$ is isomorphic to the group of
$G$-invariant divisors on $X$ modulo the subgroup of divisors of $G$-invariant
meromorphic functions.

{\sl Proof.} Let $D = \sum_{x\in X} n_xx$ be a $G$-invariant divisor. This means that,
for any $g\in G$, 
$$D = g^*(D) = \sum_{x\in X} n_xg^{-1}(x).$$
Let $L_D$ be the line bundle whose sheaf of sections is the invertible sheaf
$\calO_X(D)$ whose set of sections over an open subset $U$ is equal to $\{f\in
\bbC(X):{\tenrm div}(f)+D\ge 0\quad\hbox{after restriction to
$U$}\ \}$. The group
$G$ acts naturally on the field $\bbC(X)$ of meromorphic functions on
$X$. If $f\in \bbC(X)$ is considered as a holomorphic map $f\to \bbP^1$ then the
image 
$^gf$ of $f$ under $g\in G$ is equal to the composition $f\circ g^{-1}$. Since
$({}^gf) = g^*((f))$ we have a natural isomorphism of invertible sheaves
$\calO_X(D)\to \calO_{X}(g^*(D))$. It defines an isomorphism of line bundles 
$\phi_g:g^*(L_D)\to L_D$ which satisfies $\phi_{g'\circ g} = \phi_g\circ g^*(\phi_{g'})$.
This makes $L_D$ a $G$-linearized line bundle. If $L_D$ is equal to zero in $Pic(X)$
then 
$D = (f)$ for some $f\in \bbC(X)$ with the property $({}^gf) = (f)$ for all $g\in G$. The
ratio
$\chi_g = {}^gf/f$ is a nonzero constant, and the map $G\to \bbC^*$ defined by $\chi_g$
is a homomorphism of groups. It defines a linearization on $L_D$. It is trivial if and only if $f\in
\bbC(X)^G$. This shows that the group
${\tenrm Div}(X)^G/{\tenrm div}(\bbC(X)^G)$ of
$G$-invariant divisors  modulo principal divisors of the form $(f), f\in \bbC(X)^G,$ is mapped
isomorphically onto a subgroup of
$Pic(G;X)$. I claim that the image is the whole group. In fact, let $L$ be a $G$-linearized line bundle
and
$\phi_g:g^*(L)\to L$ be the set of isomorphisms satisfying $\phi_{gg'} = \phi_g\circ
g^*(\phi_{g'})$ which define the linearization. Then $\phi_g$ is defined by a meromorphic function $f_g$
such that $g^*(D) = D+(f_g)$. We have $f_{g'\circ g} = {}^gf_{g'}f_g$ so that
$(f_g)_{g\in G}$ is a one-cocycle of $G$ with values in $\bbC(X)^*$. By Hilbert's
theorem 90 this cocycle must be trivial. Hence we can write $f_g = {}^ga/a$ for
some $a\in \bbC(X)$. Replacing $D$ with $D' = D-(a)$ we obtain $g^*(D') = D'$ for
any $g\in G$. This shows that $L_D\cong L_{D'}$ arises from a $G$-invariant
divisor. This proves the assertion.

\medskip
\plainproclaim Proposition 2.2. There is an exact sequence of abelian groups
$$0\to Hom(G,\bbC^*)\to Pic(G;X)\to Pic(X)^G\to H^2(G,\bbC^*)\to 0.$$

{\sl Proof.} The only non-trivial assertion here is the surjectivity of the map
$$e:Pic(X)^G\to H^2(G,\bbC^*).$$ 
To prove it we need a homological interpretation
of the exact sequence. We use the following two spectral sequences with the same
limit (see {\bf [Gr]}):
$$'E_2^{pq} = H^p(G,H^q(X,\calO_X^*))\Rightarrow H^n(G;X,\calO_X^*),$$
$${}''E_2^{pq} = H^p(Y,R^q\pi_*^G(\calO_X^*))\Rightarrow H^n(G;X,\calO_X^*).$$
Here $\pi:X\to Y = X/G$ is the canonical projection and the group $H^1(G;X,\calO_X^*)$ is
isomorphic to
$Pic(G;X)$. The first spectral sequence gives the exact sequence
$$0\to Hom(G,\bbC^*)\to Pic(G;X)\to Pic(X)^G\to H^2(G,\bbC^*)\to H^2(G;X,\calO_X^*).$$
In order to show that $H^2(G;X,\calO_X^*) = 0$ we use the second spectral sequence.
We have 
$$'E_2^{20} = H^2(Y,\pi_*^G(\calO_X^*)) = H^2(Y,\calO_X^*) = 0,$$
as it follows from the exponential exact sequence $0\to \bbZ\to \calO_Y\to
\calO_Y^*\to 0$.
$$'E_2^{11} = H^1(Y,R^1\pi_*^G(\calO_X^*)) = 0$$ since
$R^1\pi_*^G(\calO_X^*)$ is concentrated at a finite set of branch points of
$\pi$.
$$'E_2^{02} = H^0(Y,R^2\pi_*^G(\calO_X^*)) = 0$$
since $R^2\pi_*^G(\calO_X^*)_y \cong H^2(G_x,\bbC^*)$, where $G_x$ is the isotropy
group of a point $x\in \pi^{-1}(y)$, and the latter group is trivial because
$G_x$ is a  cyclic group. All of this shows that $H^2(G;X,\calO_X^*) = 0$.

\medskip
\plainproclaim Corollary 2.3. Let Div$(X)^G$ be the group of $G$-invariant divisors and
$P(X)^G$ be its subgroup of principal $G$-invariant divisors. Then ${\tenrm
Div}^G(X)/P(X)^G$ is isomorphic to a subgroup of $Pic(X)^G$ and the quotient group is
isomorphic to
$H^2(G,\bbC^*)$.

\medskip
Now let us use the second spectral sequence for the map $\pi:X\to X/G$ to compute
$Pic(G;X)$ more explicitly. Let $y_1,\ldots,y_n$ be the branch points of $\pi$ and
$e_1,\ldots,e_n$ be the corresponding ramification indices. For each point $x\in
\pi^{-1}(y_i)$ the stabilizer subgroup $G_x$ is a cyclic group of order $e_i$. The
exact sequence arising from the second spectral sequence looks as follows:
$$0\to Pic(Y)\to Pic(G;X)\to \oplus_{i=1}^n \bbZ/e_i\bbZ\to 0.\eqno (2.1)$$
Here the composition of the first homomorphism with the forgetting map $e:Pic(G;X)\to
Pic(X)$ is the natural map
$\pi^*:Pic(Y)\to Pic(X)$. The second homomorphism is defined by the local isotropy
representation $\rho_{x_i}:G_{x_i}\to \bbC^*$ defined by the $G$-linearized bundle
$L$. Here we fix some $x_i$ in each fibre
$\pi^{-1}(y_i)$. Let
$D_i =
\pi^{-1}(y_i)$ considered as a reduced $G$-invariant divisor on $X$.  Let us assume that

$$Y =
\bbP^1.$$ 
Then the isomorphism classes $s_i$ of $L_i = L_{D_i}, i = 1,\ldots,n,$ 
generate
$Pic(G;X)$ and satisfy the relations $e_1s_1 =\ldots = e_ns_n$. This easily implies that
$$Pic(G;X) \cong \bbZ\oplus\bbZ/d_{1}\oplus \bbZ/(d_2/d_1)\oplus\ldots\oplus
\bbZ/(d_{n-1}/d_{n-2}),\eqno (2.2)$$ where
$$d_1 = (e_1,\ldots,e_n), d_2 = (\ldots,e_ie_j,\ldots,),\ldots, d_{n-1} =
(e_1\cdots e_{n-1},\ldots,e_2\cdots e_n).$$
To define a generator of the free part of $Pic(G;X)$ we use the Hurwitz formula for the
canonical line bundle of $X$:
$$K_X = \pi^*(K_Y)\otimes\big(\otimes_{i=1}^n L_i^{e_i-1}\big).$$ 
We know that $Pic(Y)$ is generated by the isomorphism class $\alpha$ of the line bundle
$L_{y}$ corresponding to the divisor $D = 1\cdot y$ where $y\in Y$. Then 
$\pi^*(\alpha) = N\gamma$, where $N$ is a positive integer and $\gamma$ generates
$Pic(G;X)$ modulo the torsion subgroup. Applying (2.1) we obtain that 
$$\bbZ/N\oplus Tors (Pic(G;X)) \cong \oplus_{i=1}^n \bbZ/e_i\bbZ,$$  
and hence
$$N = {e_1\cdots e_n\over (e_1\cdots e_{n-1},\ldots,e_2\cdots e_n)} =
l.c.m.(e_1,\ldots,e_n).$$ 
Now, switching to the additive notaton, 
$$e_1\ldots e_n K_X = e_1\ldots e_n\pi^*(-2\alpha)+\sum_{i=1}^n (e_i-1)e_1\ldots
e_nL_i =$$
$$=-2e_1\ldots e_nN\gamma+\sum_{i=1}^n (e_i-1){e_1\ldots
e_n\over e_i}N\gamma = e_1\ldots e_n (n-2-\sum_{i=1}^ne_i^{-1})N\gamma.$$
This implies that
$$K_X = {(n-2)e_1\ldots e_n-e_1\cdots e_{n-1}-\ldots e_2\cdots e_n\over (e_1\cdots
e_{n-1},\ldots,e_2\cdots e_n)}\gamma = l.c.m (e_1,\ldots,e_n)(n-2-\sum_{i=1}^n{1\over
e_i})\gamma.\eqno (2.3)$$ Now we are ready to compute $Pic(X(p))^G$, where $X(p)$ is
the modular curve of level
$p$ and $G = PSL(2,\bbF_p)$.  
\medskip

The following result is contained in {\bf [AR]}, Corollaries 24.3 and 24.4. However, keeping in mind
some possible applications to more general situations, we shall give it another 
proof which is based on the 
previous  discussion.

\bigskip
\plainproclaim Theorem 2.4. Assume $p \ge 5$ is prime. Let $G = PSL(2,\bbF_p)$.
Then
$$Pic(X(p))^G = Pic(SL(2,\bbF_p);X(p)) = \bbZ\lambda,$$
where 
$$\lambda^{2p-12} = K_{X(p)}$$
and 
$${\tenrm deg} \lambda = {p^2-1\over 24}.$$
Moreover, $Pic(G;X(p))$ is the subgroup of $Pic(X(p))^G$ generated by $\lambda^2$. 

{\sl Proof.} We use that the map $\pi:X(p)\to X(p)/PSL(2,\bbF_p)$ is
ramified over three points with ramification indices $2,3$ and $p$. It follows from the
previous computation that 
$Pic(G;X(p))$ is a free cyclic group generated by a
$(p-6)$-th root of the canonical class. We use Proposition 2.1 and well-known facts that
$Hom(G,\bbC^*) = \{1\}, H^2(G,\bbC^*) \cong \bbZ/2$. This gives us that $Pic(X(p))^G =
Pic(G;X(p))$ is a subgroup of index 2 in $Pic(X(p))^G$. It remains to show that 
$Pic(X(p))^G$ does not contain $2$-torsion elements. Let $L\in Pic(X)_n^G$ be a torsion element of order $n$
in
$Pic(X)^G$. Let $\mu_n$ be the constant sheaf of $n$-th roots of unity. The Kummer sequence  
$$0\to \mu_n\to \calO_X^*\to \calO_X^*\to 0$$
implies that $Pic(X)_n^G = H^1(X,\mu_n)^G$. Replacing $\calO_X^*$ with $\mu_n$ in the proof of
Proposition 2.2, we obtain that all arguments extend to this situation except that we cannot use that 
$H^2(G_x,\mu_n) = 0$. As a result we obtain an exact sequence
$$0\to Pic(G;X)_2\to Pic(X)_2^G\to H^2(G,\mu_2) \to H^2(G_{x_1},\mu_2).$$
The last homomorphism here is the restriction homomorphism for group cohomology. Here $x_1$ is a
ramification point of index 2. The exact sequence 
$$0\to \mu_2 \to \bbC^*\to \bbC^*\to 0$$
shows that $H^2(G,\mu_2) = H^2(G,\bbC^*)_2 = \bbZ/2$ (because $Hom(G,\bbC^*) =
1$) and 
$H^2(G_{x_1},\mu_2) = \bbZ/2$ (because $H^1(G_{x_1},\bbC^*) = \bbZ/2$). I claim that
the restriction homomorphism is bijective. Let $\alpha$ be the non-trivial element of
$H^2(G,\mu_2)$. It is represented by the extension 
$$1\to\mu_2\to SL(2,\bbF_p)\to
PSL(2,\bbF_p)
  \to 1.$$
 Let $\tilde
g\in SL(2,\bbZ)$ be a lift of the generator $g$ of $G_{x_1}$. Then $\tilde g^2 = -1$ 
has order 4 and therefore the exact sequence restricts to the nontrivial
extension
$$1\rightarrow\mu_2\rightarrow \mu_4\rightarrow G_{x_1}  \rightarrow 1$$
representing the nontrivial element of $H^2(G_{x_1},\mu_2)$.

The last
assertion follows from the known degree of the canonical class of a modular curve (see
{\tenbf [Sh]}).

\bigskip\noindent
{\bf Remark 2.5} The previous result implies that the group of $PSL(2,\bbF_p)$-invariant
divisors on $X(p)$ modulo principal divisors is generated by a divisor of degree
$(p^2-1)/12$. This result can be also found in ({\tenbf [AR]}, Corollary 24.3) together
with an  explicit representative of this class
$$D = \epsilon(D_2-D_3-pD_p),$$
where $\epsilon = \pm 1$ and $p = 6n+\epsilon$. Here $D_k$ denote the $G$-orbit of
points with isotropy subgroup of order $k$.

\bigskip
The tensor powers of the line bundle $\lambda$ generating the group $Pic(X(p)^G$ allows
one to embed $X(p)$ $SL(2,\bbF_p)$-equivariantly in  projective space. We state without
the proof the following result (see {\tenbf [AR]}, Corollary 24.5):

\medskip
\proclaim {Theorem 2.6} Assume $p$ is prime $\ge 5$. Denote by $V_-$ (resp. $V_+$) one
of the two irreducible representations of $SL(2,\bbF_p)$ of dimension ${p-1\over2}$
(resp. ${p+1\over 2}$). Then
\roster
\item  a base-point free linear subsystem of
$|\lambda^{(p-3)/2}|$ maps $X(p)$ in $\bbP(V_-) = \bbP^{(p-3)/2}$ onto a curve of
degree $(p-3)(p^2-1)/48$;
\item  a base-point-free linear subsystem of $|\lambda^{(p-1)/2}|$ maps
$X(p)$ in
$\bbP(V_+) = \bbP^{(p-1)/2}$ onto a curve of degree $(p-1)(p^2-1)/48$.
\endroster
\endproclaim
\medskip
It is conjectured that the linear systems embedding $X(p)$ in $\bbP(V_-)$ and in $\bbP(V_+)$ are complete 
(see {\tenbf [AR]}, p.106). This is known to be true only for $p = 7$. 

\bigskip\noindent
{\bf Remark 2.7} As was shown in {\tenbf [AR]}, Corollary 24.5, the image of $X(p)$ in
$\bbP(V_+)$ (resp. 
$\bbP(V_-)$ ) described in the previous theorem is the $z$-curve (resp.
$A$-curve) of Klein. From the modern point of view these embeddings can be
described as follows. Recall that $X(p)$ is a compactification of the moduli
space of isomorphism classes of pairs $(E,\phi)$, where $E$ is an elliptic curve
and $(e_1,e_2)$ is a basis of its group $E_p \cong (\bbF_p)^2$ of $p$-torsion
points. Let $O$ be the origin of $E$. There is a special basis
$(X_0,\ldots,X_p)$ in the  space $\Gamma(E,\calO_E(pO))$ which defines a map
$f:E\to \bbP^{p-1}, x\to (X_0(x):\ldots:X_p(x))$ satisfying the following
properties:
\roster
\item $(X_0(x-e_1):\ldots:X_p(x-e_1)) = (X_1(x):X_2(x):\ldots:X_p(x):X_0(x))$;
\item $(X_0(x+e_2):\ldots:X_p(x+e_2)) =
(X_0(x):\zeta X_1(x):\ldots:\zeta^pX_p(x))$, where $\zeta = e^{2\pi i/p}$;
\item $(X_0(-x):\ldots:X_p(-x)) = (X_p(-x):X_{p-1}(x):\ldots:X_1(x):X_0(x))$.
\endroster
Let $\pi:\calX(p)\to X(p)$ be the universal family of elliptic curves $(E,e_1,e_2)$ (its
fibres over cusps are certain degenerate curves, $p$-gons of lines). The $p$-torsion
points of the fibres determine $p^2$ sections of the elliptic surface $\calX(p)$. The
functions
$X_m$ define a morphism $\calX(p)\to \bbP^{p-1}$ whose restriction to the fibre
$(E,e_1,e_2)$ is equal to the map $f$. The functions $z_m = X_m-X_{p-m}, m = 1,\ldots,
{p-1\over 2},$ define a projection of the image of $\calX(p)$ in $\bbP^{(p-3)/2}(\bbC)$.
The image of the section of $\pi$ defined by the $0$-point is the $z$-curve of Klein. On
the other hand if we consider the functions $y_m = X_m+X_{p-m}$ we get the projection of
$\calX(p)$ in
$\bbP^{(p-1)/2}(\bbC)$. The center of this projection contains the $0$-section so
that the restriction of the projection map to each curve $f(E)$ is the projection
from the origin $f(0)$ and still defines the image of $f(0)$. The set of these images
forms the $A$-curve of Klein. There are certain modular forms $A_0,\ldots,A_{{p+1\over
2}}$ defined on the upper-half plane which define the map from $X(p)$ to the $A$-curve.

It is proven in {\tenbf [Ve]} that the $z$-curve is always nonsingular.

\bigskip\noindent
{\bf Example 2.8.} Let $p = 7$. Then deg $\lambda = 2$, $K_X = \lambda^2$ is of
degree $4$ so that $X(7)$ is a curve of genus $3$. The divisor class $\lambda$ is an even
theta-characteristic on $X(7)$. It is the unique theta characteristic invariant with
respect to the group of automorphisms $G = PSL(2,\bbF_7)$ of $X(7)$ (see other proofs 
of this fact in {\tenbf
[Bu], [DK]}). The
$z$-curve is a canonical model of
$X(7)$, a plane quartic. This is the famous Klein's quartic with $168$ automorphisms. In
an appropriate coordinate system it is given by the equation
$$x_0^3x_1+x_1^3x_2+x_2^3x_0 = 0.$$
The $A$-curve is a space sextic with equations
$$t^2x+\sqrt{2} ty^3+2yz^2 = t^2y+\sqrt{2} tz^2+2zx^2 = t^2z+\sqrt{2} tx^2+2xy^2 = 2\sqrt{2}
xyz-t^3 = 0$$
(see {\tenbf [E1]}, p. 163).

It is well-known (see, for example, {\tenbf [Tj]}) that a theta characteristic $\theta$
with $H^0(X,\theta) = 0$ on a plane nonsingular curve $X$ of degree $n$ with equation
$F = 0$ gives rise to a representation of $F$  as the determinant of a symmetric $n\times
n$ matrix with linear forms as its entries. In other words, $\theta$ defines a net of
quadrics in $\bbP^{n-1}$ and $X$ parametrizes singular quadrics from the net. The pair
$(X,\theta)$ is called the Hesse invariant of the net. It follows from Table 2 in Appendix
1 that
$S^2(V_+)$ contains
$V_-$ as a direct summand (as representations of $SL(2,\bbF_7))$. This defines a $SL(2,\bbZ)$-invariant net
of quadrics in
$\bbP(V_+)$ with the Hessian invariant $(X,\lambda)$. The corresponding
representation of $X$ as the determinant of a symmetric matrix of linear forms is known
since Klein (see {\bf [E1]}, p. 161):
$$x_0^3x_1+x_1^3x_2+x_2^3x_0  = \det\pmatrix -x_0&0&0&-x_1\\
0&x_1&0&-x_2\\
0&0&x_2&-x_0\\
-x_1&-x_2&-x_0&0\endpmatrix.\eqno (2.4)$$

\bigskip\noindent
{\bf Example 2.9.} Let $p = 11$. Then deg $\lambda = 5$, $K_X = \lambda^{10}$ is of
degree $50$ so that $X(11)$ is a curve of genus $26$. The $z$-curve is a curve of degree
$20$ in $\bbP^4$. According to F. Klein {\tenbf [Kl]} it is equal to the 
locus 
$$\{(v:w:x:y:z)\in \bbP^4: rk\pmatrix w&v&0&0&z\\
v&x&w&0&0\\
0&w&y&x&0\\
0&0&x&z&y\\
z&0&0&y&v\endpmatrix \le 3\}.$$
The matrix in above is the Hessian matrix of a $G$-invariant cubic hypersurface $W$
given by the equation
$$v^2w+w^2x+x^2y+y^2z+z^2v = 0.$$
This hypersurface has the group of automorphisms isomorphic to $PSL(2,\bbF_{11}).$ 
The $A$-curve is a curve of degree 25 in $\bbP^5$. It is the curve of
singularities of a unique quartic ruled hypersurface in $\bbP(V_+)$ 
(see Appendix V in {\tenbf [AR]} which contains the results of the first author). We
refer for these and other beautiful facts about the geometry of $X(11)$ to {\tenbf
[AR]} and {\tenbf [E2]}.

\bigskip
\noindent
{\bf 3. Rank 2 bundles}. We shall use the following result of S. Ramanan which is a special case of
Proposition 24.6 from {\bf [AR]}:

\plainproclaim Theorem 3.1. Let $G$ be a finite subgroup 
of the group of automorphisms of a curve $X$ and $E$ be a $G$-linearized
rank $r$ vector bundle over $X$. Then there exists a flag
$$0\subset E_1\subset E_2\subset\ldots \subset E_{r-1}\subset E$$
of $G$-invariant subbundles, where each $E_i$ is of rank $i$ and all inclusions are
$G$-equivariant.

\medskip
Because of the importance of this result for the sequel we shall
sketch a proof. Choose an ample $G$-linearized line bundle $L$ with trivial isotropy representations. This
is always possible by taking products and powers of $g$-translates of an ample line bundle. Then we apply
the Lefschetz Fixed Point Formula for coherent sheaves
$$tr(g|H^0(X,E\otimes L^n))-tr(g|H^1(X,E\otimes L^n)) = \sum_{g(x) = x}{tr(g|L_x^n\otimes
E_x)\over 1-dg_x}.$$ Since $g$ acts identically on $L_x$, the right-hand side is
independent of $n$. Taking $n$ sufficiently large, we get rid of $H^1$. Since the dimension
of $H^0$ will grow with $n$ and the trace of $g\ne 1$ on
$H^0$ does not change with $n$ we easily obtain that $H^0$ contains the trivial irreducible representation
for large
$n$. This implies that there exists a $G$-invariant section of $E\otimes L^{n}$. It gives a $G$-invariant
embedding of  $L^{-n}$ in $E$. It generates a $G$-invariant subbundle of $E$. Now we take the quotient
and apply induction on the rank.

\medskip
\plainproclaim Corollary 3.2. Let $E$ be a $PSL(2,\bbF_p)$-invariant rank $r$ vector bundle
over $X(p)$. Then there is a $PSL(2,\bbF_p)$-equivariant flag
$$0\subset E_1\subset E_2\subset\ldots \subset E_{r-1}\subset E \eqno (3.1)$$
of $PSL(2,\bbF_p)$-invariant rank $i$ vector bundles $E_i$. 

{\sl Proof.} To see that the previous theorem applies, we use Corollary 1.2 that shows
that $E$ admits a $SL(2,\bbF_p)$-linearization. If the center $C =\{\pm 1\}$ of
$SL(2,\bbF_p)$ does not act identically on $E$, hence acts as $-1$ on each fibre, we
replace $E$ by $E' = E\otimes
\lambda$. Since $\lambda$ is $SL(2,\bbF_p)$-linearized but not
$PSL(2,\bbF_p)$-linearized the center $C$ acts as $-1$ on $\lambda$. So $E'$ is
$PSL(2,\bbF_p)$-linearized and the theorem applies. 

\medskip
We shall call a flag (3.1) a {\it Ramanan flag} of $E$.
Let $L_i = E_i/E_{i-1}, i = 1,\ldots,r$, where $E_0 = 0, E_r = E$, be the factors of a
Ramanan flag of $E$. We know that each $L_i$ is equal to $\lambda^{a_i}$ for some
integer $a_i$. The shall call the sequence $(a_1,\ldots,a_r)$ a {\it sequence of
exponents} of
$E$. Note that it may not be defined uniquely.

\medskip
\plainproclaim Proposition 3.3. Let $(a_1,\ldots,a_r)$ be a sequence of exponents of a
$PSL(2,\bbF_p)$-invariant stable rank $r$ bundle over $X(p)$. Let $a =
a_1+\ldots+a_r$. Then, for any $s < r$,
 $$a_1+\ldots+a_s < {s\over r}a.$$

{\sl Proof.} This follows immediately from the definition of stability. 

\medskip
In the case $r = 2$ we will be able to say 
more about sequences of exponents of a rank 2 bundle (see Corollary 4.3) but now let us note the
following result (see {\bf [AR]},  Lemma 24.6):

\medskip
\plainproclaim Proposition 3.4. Assume $r = 2$ and let $(a_1,a_2)$ be a sequence of exponents of a $G$-stable
bundle $E$. Then $a_1+a_2$ is even.

\medskip
This follows from the fact that any $G$-invariant extension
$$0\to \lambda^{a_1}\to E\to \lambda^{a_2}\to 0$$
has obstruction class for splitting in 
$$H^1(X,\lambda^{a_1-a_2})^{SL(2,\bbF_p)} =
(H^0(X,\lambda^{a_2-a_1+2p-12})^*)^{SL(2,\bbF_p)}.$$ 
Since $-1$ acts as $-1$ in $H^0(X,\lambda^{{\tenrm
odd}})^*$ the latter space is trivial.

\medskip
\plainproclaim Corollary 3.5. Each $PSL(2,\bbF_p)$-invariant stable bundle of rank 2 over
$X(p)$ has determinant isomorphic to $\lambda^a$, where $a$ is even.

\medskip
By tensoring $E$ with $E\otimes \lambda^{-a/2}$ we may assume now that $\det E$ is
trivial. This allows us to invoke some results from topology. Recall the
following

\medskip
\plainproclaim Theorem 3.6 ({\tenbf [NS]}). Let $E$ be a degree $0$ semi-stable vector
bundle on a compact Riemann surface $X$. Then
there exists a unitary representation
$$\rho :\pi_1(X)\to U(r)$$
such that $E$ is isomorphic to the vector bundle $H\times \bbC^r/\pi_1 \to X =
H/\pi_1(X)$, where $H$ is the universal cover of $X$ with the natural action of
$\pi_1(X)$ on it, and the fundamental group $\pi_1$ acts on the product by the formula
$\gamma:(z,v)\to (\gamma\cdot z,\rho(\gamma)\cdot v)$. This construction defines a
bijective correspondence between the set of isomorphism classes of semi-stable rank $r$
bundles of degree zero and the set of unitary representations of $\pi_1(X)$ of
dimension $r$ up to conjugation by a unitary transformation of $\bbC^r$. In this 
correspondence stable bundles correspond to
irreducible representations and semi-stable bundles with trivial determinant correspond
to representations $\rho: \pi_1(X)\to SU(r)$.

\medskip
We shall apply this theorem to our situation. First we need the following:

\medskip\noindent
{\bf Definition.} A $G$-linearized vector bundle is called {\it $G$-stable} if for any
$G$-linearized subbundle $F$ of $E$ one has $\mu(F) < \mu(E)$.

\medskip
Notice that a $G$-stable bundle is always semi-stable since it is known that an unstable
bundle has always a unique, hence $G$-invariant, maximal destabilizing subbundle (see
{\bf [Se]}). Obviously a stable bundle is $G$-stable. However, a semi-stable
bundle could be also
$G$-stable. For example, the trivial bundle defined by an irreducible representation of
$G$ is $G$-stable but not stable.

\medskip
\plainproclaim Theorem 3.7.  Let $E$ be a $G$-linearized semi-stable vector bundle of degree
$0$ given by a unitary representation $\rho :\pi_1(X)\to U(r)$. Let $\Pi$ be the group of
automorphisms of the universal cover of
$X$ generated by lifts of elements of $G$. Then $\rho$ can be extended to a unitary
representation $\tilde \rho:\Pi\to U(r)$. Moreover there is a bijective correspondence
between the set of isomorphism classes of $G$-linearized semi-stable rank $r$ bundles of
degree $0$ and the set of unitary representations $\Pi\to U(r)$ up to conjugation by a
unitary transformation of $\bbC^r$. In this correspondence irreducible
representations correspond to $G$-stable bundles.

{\sl Proof.} The unitary representation $\rho:\pi_1(X)\to U(r)$ arises from the
holonomy representation of a flat unitary connection on $E$ compatible with a holomorphic
structure on $E$. As soon as we fix a $G$-invariant Hermitian metric on $E$, such a
connection is defined uniquely up to isomorphism. This shows that this connection is
$G$-invariant. It is also
$G$-linearized in the sense that the corresponding local system $\calV$ of horizontal
sections is $G$-linearized. Since $E = \calO_X\otimes_\bbC {\calV}$ it is easy to see
that the obstruction for the linearization of ${\calV}$ is determined by the same
element
$c(E)$ as the obstruction for the linearization of
$E$. By assumption the latter is trivial. Now the group $\Pi$ acts on the pull-back
$\tilde {\calV}$ on $H$ by $(z,v)\to (pz, \tilde\rho(p)(v))$ for some linear
representation $\tilde\rho:\Pi\to U(r)$ extending the representation $\rho$.

Conversely, given a representation $\rho:\Pi\to U(r)$, restricting it to the subgroup
$\Gamma = \pi_1(X)$ we get a unitary representation $\rho$ of $\pi_1(X)$ and hence a
semi-stable bundle $E = H\times \bbC^r/\pi_1(X)$ on $X$. The factor group $G
=\Pi/\pi_1(X)$ acts naturally on $E$ preserving the structure of a vector bundle. This
action defines a $G$-linearization. 

\medskip
Let $\Pi$ be as above. It is given by an extension of groups
$$1\to \Gamma \to \Pi \to G\to 1,$$
where $\Gamma\cong \pi_1(X)$. Assume that $X/G = H/\Pi = \bbP^1$. Then 
$\Pi$, as an abstract group, is given by the genetic code:
$$\Pi = <\gamma_1,\ldots,\gamma_n|g_1^{e_1} = \ldots = \gamma_n^{e_n} = \gamma_1\cdots
\gamma_n = 1>,$$ where $\sum e_i^{-1} < n-2$. As a group of transformations of $H$, $\Pi$
is isomorphic to a discrete subgroup of $PSL(2,\bbR)$ which acts on $H$ as
a subgroup generated by even products of reflections in sides of a geodesic $n$-gon
with  angles 
$\pi/e_i$. Such a subgroup of $PSL(2,\bbR)$ is called a {\it Dyck group} (or a {\it
triangle group} if $n = 3$) of signature $(e_1,\ldots,e_n)$. Conversely, let
$\Pi$ be an abstract group as above and
$\Gamma$ be a normal torsion-free subgroup of finite index. Choose an isomorphism
from
$\Pi$ to
a Dyck group (if $n = 3$ it is defined uniquely, up to a conjugation). Let us identify
$\Pi$ with its image. The group $\Gamma$ acts freely on $H$ and the quotient $H/\Gamma$
is a compact Riemann surface $X_\Gamma$ with
$\pi_1(X) \cong \Gamma$. The factor group $G = \Pi/\Gamma$ acts on $X_\Gamma$ by
holomorphic automorphisms. The projection
$\pi:X_\Gamma\to X_\Gamma/G = 
\bbP^1$ ramifies over $n$ points with ramification indices $e_1,\ldots,e_n$. 

The
modular curve
$X(p)$ is a  special case of this construction. One takes $(e_1,\ldots,e_n) =$
$(2,3,p)$
and $\Gamma =\pi_1(X(p))$.

Let $\tilde{\Pi}$ be the
group with the genetic code 
$$\tilde{\Pi} = <\tilde \gamma_1,\ldots,\tilde \gamma_n,h|h\ \  \hbox{central}, \tilde \gamma_1^{e_1} =
\ldots =
\tilde \gamma_n^{e_n} = \tilde \gamma_1\ldots \tilde \gamma_n = h>.
$$
 It is a central extension of $\Pi$ with infinite cyclic center
generated by $h$:
$$1\to (h)\to \tilde{\Pi}\to \Pi \to 1.$$

We shall assume that $\Pi$ is perfect, i.e. coincides with its commutator
$\Pi' = [\Pi,\Pi]$. This happens if and only if the $e_i$'s are pairwise 
prime.  In this case the commutator group
$\tilde \Pi' = [\tilde{\Pi},\tilde{\Pi}]$ is a universal central extension of $\Pi$:
$$1\to (t)\to [\tilde{\Pi},\tilde \Pi]\to \Pi \to 1.\eqno (3.3)$$
 Its center $(t)$ is a subgroup of
$(h)$ generated by $h^s$, where 
$$s = e_1\cdots e_n(n-2-\sum e_i^{-1}).$$
The genetic code of $\tilde{\Pi}'$ is
$$\tilde {\Pi}' = <\tilde g_1,\ldots,\tilde g_n,t|t\quad\hbox{central}, \tilde g_1^{e_1} = t^{b_1},
\ldots, \tilde g_n^{e_n} = t^{b_n}, \tilde g_1\cdots \tilde g_n = t^b>,$$ where 
$$sb_i \equiv 1\ {\tenrm mod}\  e_i,\quad 0\le b_i < e_i,$$
$$b = {1\over e_1\cdots e_n}+\sum_{i=1}^n{b_i\over e_i}.$$
All of this is well-known in 3-dimensional topology (see {\bf [FS]}) and the theory of singularities of
algebraic surfaces (see {\bf [Do]}). We have
$$\tilde \Pi' = \pi_1(\Sigma(e_1,\ldots,e_n)),$$
where $\Sigma(e_1,\ldots,e_n)$ is a Seifert-fibred 3-dimensional homology sphere given
explicitly as the intersection of a sphere $S^{2n-1}$ with center at the origin in
$\bbC^{n}$ and the algebraic surface given by the equations
$$z_1^{e_1}+z_i^{e_i}+z_n^{e_n} = 0,\quad i= 2, \ldots, n-1.$$
The group $\Pi$ is the fundamental group of the link of a canonical Gorenstein
singularity admitting a good $\bbC^*$-action (see {\bf [Do]}).

\bigskip
\plainproclaim Corollary 3.8. Keep the notation of Theorem 3.7. Assume that the group $\Pi$ is
perfect. Let
$1\to H^2(G,\bbC^*)\to \tilde G\to G\to 1$ be the universal central extension of $G$
and let $d = |H^2(G,\bbC^*)|$. Then $\tilde\Pi'$ is mapped surjectively on $\tilde
G$ and there is a bijective correspondence between irreducible unitary representations
$\rho:\tilde
\Pi'\to SU(r)$ with
$\rho(h)^d = 1$ which restrict to an irreducible representation of $\tilde\Gamma =
Ker(\tilde \Pi'\to \tilde G)$ and
$G$-invariant stable rank $r$ bundles over $X_\Gamma$ with trivial determinant.

{\sl Proof.} Let $\Pi \cong  F/R$ where $F$ is a free group. It is known that the
universal central extension $\tilde \Pi'$ is isomorphic to $[F,F]/[R,F]$ ( see for
example {\bf [Mi]}). Since $G
\cong F/R'$, where $R\subset R'$ we obtain that $\tilde G \cong [F,F]/[R',F]$ and there is a
surjective homomorphism  
$\tilde
\Pi'\to \tilde G$. Let $\tilde \Gamma$ be the kernel of this homomorphism.  We  have a central
extension for
$\tilde\Gamma$:
$$1\to (t^{d'}) \to \tilde \Gamma \to \Gamma\to 1,$$ 
where $d'|d$. Given an irreducible representation 
$\rho:\tilde
\Pi'\to SU(r)$ with $\rho(h)^d = 1$ we define $\beta:\Gamma\to SU(r)$ by first
restricting
$\rho$ to $\tilde\Gamma$ and then factoring it through the quotient $\tilde
\Gamma/(t^d)
\cong
\Gamma$. Since, by the assumption, the restriction of $\rho$ to $\tilde \Gamma$ is irreducible,
$\beta$ is irreducible. This defines a stable rank
$r$ bundle
$E$ on
$X_\Gamma$ with trivial determinant. Since $\tilde\Gamma$ is normal in $\tilde \Pi'$, the
group $\tilde G$ acts on $E = H\times \bbC^r/\tilde \Gamma$ and makes $E$ a
$G$-invariant bundle. Conversely, by Corollary 1.2, any $G$-invariant stable bundle 
$E$ with trivial determinant admits a $\tilde G$-linearization. This linearization
defines a $\tilde G$-linearization on the local coefficient system $\calV$ defined by
the flat unitary connection on $E$. The group $\tilde \Gamma$ acts on the pull-back 
$\bar{\calV} = \calV \times_XH$
of $\calV $ on $H$ via the action of $\Gamma$ on $H$ and the trivial action on
${\calV} $. Since $\bar{\calV}$ trivializes on $H$, it coincides with the universal
covering of
${\calV} $. Thus the action of $\tilde G$ lifts to an action of $\bar{\calV}$ and
defines an action of $\Pi'$ on $\bar{\calV} = H\times \bbC^r$. This defines a unitary
representation of the  group
$\tilde
\Pi'$ in $SU(r)$.
 
\medskip\noindent
{\bf Remark 3.9} If $r = 2$ we do not need the assumption that the restriction $\beta$ of
$\rho:\tilde \Pi'\to SU(2)$ to $\tilde \Gamma$ is irreducible.  Assume it is not. Then
it decomposes into sum of one-dimensional unitary representations $V_1+V_2$. By
Theorem 3.6 they define two line bundles on $X_\Gamma$ with determinants of degree 0.
These bundles are invariant with respect to the
group $G = \Pi/\Gamma$. However, the computations from the first section show that
$Pic(X_\Gamma)$ is generated by an element of positive degree. This shows that $\beta$
must be the trivial representation. Thus $\rho$ factors through a representation $\bar
\rho:G\to SU(2)$. However, it is easy to see using the classification of finite
subgroups of $SU(2)$ that 
$G$ does not admit non-trivial 2-dimensional unitary representations. This gives us that
$\rho$ is the trivial representation which contradicts the assumption that $\rho$ is
irreducible.

\medskip
\plainproclaim Corollary 3.10. Let $(a_1,a_2)$ be a sequence of exponents of a stable
$G$-invariant rank 2 vector bundle over $X(p)$. Then $a_1$ and $a_2$ are odd
numbers.

{\sl Proof.} We already know from Proposition 3.4 that $a_1+a_2$ is even. So, it suffices to
show that $a_1$ and $a_2$ cannot be both even. Assume they are even. Then the
extension
$$0\to \lambda^{a_1}\to E\to \lambda^{a_2}\to 0$$
shows that the center of $SL(2,\bbF_p)$ acts identically on $E$. Hence the bundle $E$
admits a
$PSL(2,\bbF_p)$-linearization. By Theorem 3.7, this defines a unitary representation
$\rho:\Pi\to SU(2)$. Let $g_1$ be a generator of $\Pi$ of order $2$. Then $\rho(g_1)^2 =
1$ and hence
$\rho(g_1) =\pm 1$. This implies that $\rho(g_2)^3 = \rho(g_3)^p = 1$ and $\rho(g_2g_3)
= \rho(g_1) =\pm 1$. This gives that $\rho(g_i) = 1, i = 1,2,3$, i.e.
$\rho$ is trivial. This contradicts the assumption that $E$ is stable.

\medskip
\plainproclaim Theorem 3.11. Let $p = 6n\pm 1$. Then there exist exactly $2n$
non-isomorphic rank 2 stable $G$-invariant vector bundles over $X(p)$ with trivial
determinant.

{\sl Proof.} This is an immediate corollary of Theorem 3.7, Remark 3.9, and the known
computation of the number of irreducible unitary representation of the fundamental group of 
the Brieskorn sphere $\Sigma(e_1,e_2,e_3)$ (see {\tenbf [FS]}).

\bigskip\noindent
{\bf 4. Isotropy representation.}
Let us return to the general situation of a finite group $G$ acting on a compact
Riemann surface $X$. Let $\calC(X;G)$ be the set of pairs
$(C,g)$, where $C$ is a connected component of the fixed
locus of $g\in G$, modulo
the natural action of $G$ on these pairs by
$$g'\cdot (C,g) = (g'(C),g'gg'{}^{-1}).$$ 
Since $X$ is a curve, $C$ is either a single point or
the whole $X$. In the latter case the element $g$
belongs to the kernel $A$ of the action of $G$ on $X$. Let $\check {\calC}(X;G)$
denote the set of complex valued functions on the set $\calC(X;G)$. One defines
the {\it isotropy representation} map:
$$ \rho:SU_X(G;r)\to \check {\calC}(X;G).$$
For each $E\in SU_X(G;r)$ defined by isomorphism $\phi_g:g^*(E)\to E$ and $(C,g)\in
{\calC}(X;G)$ we have
$$\rho(E)(C,g) = Trace(\phi_{g,x}:g^*(E)_x =
E_x\to E_x),$$
where $x\in C$.

Consider
the quotient 
$$Y = X/G$$
and let $p:X\to Y$
be the natural orbit map. There is a finite set of $G$-orbits
in $X$ with non-trivial isotropy subgroup. They correspond to
the set $S$ of points $y_1,\ldots,y_n$ in $Y$ such that $p$ is
ramified over any point $x\in p^{-1}(S)$. For any $g\in G$, we
have
$gG_xg^{-1} = G_{g\cdot x}.$ 
In each fibre $p^{-1}(y_i)$ pick a point $x_i$ and denote the
corresponding isotropy subgroup $G_{x_i}$ by $G_i$. This is an
extension of a cyclic subgroup $\bar G_i$ of $\bar G$ of order
$e_i$ and the group $A$. It is clear that each
$(C,g)\in
{\calC}(X;G)$ can be represented by a pair
$(x_i,g_i), g_i\in G_i,$ or by a pair $(X,a), a\in A$. Assume $A$ is central in $G$ (as it
will be in our case). Then this representation is unique since $G_{x_i}$ is cyclic.
Thus we have
$$|{\calC}(X;G)| = \sum_{i=1}^n(e_i-1)+|A|.$$
The representatives of ${\calC}(X;G)$ can be chosen as follows:
$$(x_1,g_1),\ldots, (x_1,g_1^{e_1-1}),\ldots,(x_n,g_n),\ldots, (x_n,g_n^{e_n-1}),
(X,a), a\in A,$$ 
where each $g_i\in G$ is a representative of a generator of $(G/A)_{x_i}$. 
\medskip
Now we place
ourselves in the situation discussed in the previous section and assume that $X
= X_\Gamma$, where $\Gamma$ is a torsion free normal subgroup of a Dyck group $\Pi$ of
signature $(e_1,\ldots,e_n)$. Also we assume that $\Pi$ is perfect. In this case we can
choose representatives of ${\calC}(X_\Gamma;G)$ taking for $g_i$ the images of the
standard generators $\tilde\gamma_i$ of $\tilde\Pi'$.

\medskip
\plainproclaim  Theorem 4.1. Keep the notation of Corollary 3.8. Let $E$ be a stable vector
bundle on $X_\Gamma$ arising from a unitary representation $\rho:\tilde\Pi'\to SU(r)$.
Then, for any integer $k$,
$$Trace(\rho(\tilde g_i^k)) = Trace((\phi_{g_i^k})_{x_i}:E_{x_i}\to E_{x_i}), \quad i =
1,\ldots,n.$$

{\sl Proof}. This follows easily from the construction of $E$ by means of a unitary
representation of $\tilde\Pi'$.

\medskip
In the case $r = 2$ an algorithm for computations of the traces Trace$(\rho(\tilde
g_i^k))$ of a unitary representation $\rho$ of the group
$\tilde
\Pi'$ of signature $(e_1,e_2,e_3)$ is described in {\tenbf [FS]}. I am thankful to N.
Saveliev for realizing this algorithm for me in the case $(e_1,e_2,e_3) = (2,3,p)$.

\medskip
\plainproclaim  Theorem 4.2. Let $(e_1,e_2,e_3) = (2,3,p)$ and $\rho:\tilde\Pi'\to SU(2)$ be
an irreducible representation. Write $p = 6n+\epsilon$, where $\epsilon = \pm 1$. Then 
$$[Trace(\rho(\tilde
\gamma_1)),Trace(\rho(\tilde \gamma_2)),Trace(\rho(\tilde \gamma_3))] =
[0,\epsilon,2\cos({\pi k\over p})],$$
where $k$ is an odd integer between $n+1$ and $5n$ if $\epsilon = 1$ and an even integer
between $n$ and $5n-1$ if $\epsilon = -1$.

\medskip
Recall that we have $2n$ unitary irreducible
representations $\rho:\tilde \Pi'\to SU(2)$ and this agrees with the number of all
possible triples of the characters.

The conjugation classes of the unitary matrices $\rho(\tilde\gamma_i)$ are represented
accordingly by
$$\bigg[\pmatrix i&0\\
0&-i\endpmatrix,\pmatrix e^{(3-\epsilon)\pi i/6}&0\\
0&e^{-(3-\epsilon)\pi i/6}\endpmatrix,\pmatrix e^{k\pi i/p}&0\\
0&e^{-k\pi i/p}\endpmatrix\bigg].$$ 
So raising the corresponding matrices in powers and computing the traces, we get the
expression for the traces of powers of the generators $\tilde\gamma_i$ and of the
central element $t$.

\medskip
\plainproclaim  Corollary 4.3. Let $E$ be a  stable rank 2 $G$-invariant
bundle on $X(p)$ with trivial determinant. Let $[0,\epsilon,2\cos({\pi k\over p})]$
define its isotropy representations and let $(a,-a)$ be a sequence of
exponents of $E$. Then $a$ is an odd negative integer, and 
$$an  \equiv \epsilon k\quad\hbox{mod $2p$},\quad
a\equiv  1\quad\hbox{mod $3$}.$$

{\sl Proof.} Let
$$0\to \lambda^{a}\to E\to \lambda^{-a}\to 0$$ 
be the extension defined by the sequence of exponents $(a,-a)$. Clearly the isotropy
representation of $E$ is determined in terms of the isotropy representation of 
$\lambda^a$. We know that $\lambda^{2p-12} = K_{X(p)}$. The isotropy representation of the
cotangent line bundle $K_X$ is easy to find. Any
generator $g_i$ of the isotropy group $G_{x_i}$ acts as a primitive $e_i$-th root of unity. Let
us take it to be $e^{2\pi i/e_i}$. Then the isotropy representation of $\lambda$ at
$(x_i,g_i)$ is given by some $2e_i$-th root of unity $e^{s_i\pi i/e_i}, 0\le s_i < 2e_i,$
which satisfies
$(p-6)s_i \equiv 1 \ {\tenrm mod}\  e_i$. We easily find
$$s_i = \cases 1&\text{if $i =1$}\\
{3-\epsilon\over 2}&\text{if $i = 2$}\\
n&\text{if $i = 3, \epsilon = 1$}\\
p-n&\text{if $i = 3, \epsilon = -1$}\endcases.$$
This shows that $[2\cos(as_1\pi/2),2\cos(as_2\pi/3),2\cos(as_3\pi/p)] =
[0,\epsilon,2\cos(k\pi/p)].$ We check that the first  entries coincide automatically.
To satisfy the equality for the other two we need the congruences  stated
in the assertion.

\bigskip\noindent
{\bf 5. The Adler-Ramanan-Klein bundle.}  It
was introduced by A. Adler and S. Ramanan ({\tenbf [AR]}, \S 24). It arises from an
interpretation of Klein's quartic equations defining the
$z$-curve $X(p)$ (see {\tenbf [KF]}, p.268). We refer to {\tenbf [AR]} and
{\tenbf [Ve]} for a modern treatment of these equations.  We shall prove
that this bundle is stable when
$p = 7$ and find its sequence of exponents.

Recall that  $SL(2,\bbF_p)$ has two non-isomorphic
irreducible representations, dual to each other, of dimension
${p-1\over 2}$ and two non-isomorphic representations, dual to each other, of dimension
${p+1\over 2}$. The $z$ curve $X(p)$ lies in the projectivization of an irreducible
representation of dimension ${p-1\over 2}$. We denote it by $V_-$. The $A$-curve $X(p)$
lies in the projectivization of an irreducible representation of dimension
${p+1\over2}$. We denote it by $V_+$. It is conjectured that
$V_-^* = H^0(X(p),\lambda^{p-3/2})$ and $V_+^* = 
H^0(X(p),\lambda^{p+1/2})$ (see {\bf [AR]}, p.106). 

\medskip
\proclaim{Theorem 5.1 (A.Adler-S.Ramanan, F. Klein)} There is an isomorphism of
representations of
$SL(2,\bbF_p)$:
$$\tau:S^2(V_-) \cong \Lambda^2(V_+).$$
Let $\nu_2:\bbP(V_-)\to \bbP(S^2(V_-))$ be the Veronese embedding given by the complete
linear system of quadrics in $\bbP(V_-)$. Identifying $\bbP(S^2(V_-))$ with
$\bbP(\Lambda^2(V_+))$ by means of $\tau$, we have
$$X(p) = \nu_2^{-1}(G(2,V_+)),$$
where $G(2,V_+)$ is the Grassmann variety of 2-dimensional linear subspaces in
$V_+$.
\endproclaim
\medskip\noindent
{\bf Definition.} The {\it Adler-Ramanan-Klein bundle} (the ARK bundle for brevity) over $X(p)$
is the inverse image of the tautological rank 2 bundle over $G(2,V_+)$ under the map
\noindent
$\nu_2:X(p)\to G(2,V_+)$.

\medskip
\plainproclaim  Theorem 5.2. The determinant of the ARK bundle
$E$ is equal to $\lambda^{3-p}$. It is stable  provided  the
following condition is satisfied:
$$\hbox{$H^0(X(p),\lambda^a)$ does not contain $V_+^*$ as in
irreducible summand if $a \le {p-3\over 2}$} \eqno (*)$$

{\sl Proof.} Since the dual bundle $E^*$ embeds $X(p)$ in $G(2,V_+)$ and the
corresponding Pl\"ucker embedding of $X(p)$ is given by quadrics we obtain that the
determinant of $E^*$ is equal to $\lambda^{p-3}$. Assume
$E$ is not stable. Then $E^*$ is unstable too and contains a destabilizing subbundle
of degree
$\ge {p-3\over 2}{\tenrm deg} \lambda$. It is known that one can always choose a
unique maximal stabilizing subbundle. This implies that $E^*$ contains a $G$-invariant
subbundle isomorphic to $\lambda^a$ with $a \ge {p-3\over 2}.$ Then $E^*$ has a
quotient of the form
$\lambda^a$, where $a \le {p-3\over 2}.$ Since $E^*$ defines an embbeding in $G(2,V_+)$
it is spanned by the subspace $V_+^*$ of its space of global sections. This shows that
$\lambda^a$ is spanned by $V_+^*$ too. This implies that there is a
$SL(2,\bbF_p)$-invariant non-trivial linear map $V_+^*\to H^0(X(p),\lambda^a)$. This
contradicts the assumption of the theorem.

\medskip\noindent
{\bf Remark 5.3} It is reasonable to conjecture that condition (*) is always satisfied.
In fact, together with Adler and Ramanan, I believe that $H^0(X(p),\lambda^a) = 0$ for 
$a < {p-3\over 2}$ and $H^0(X(p),\lambda^{{p-3\over 2}})\cong V_-$
(see the ``WYSIWYG'' Hypothesis in {\tenbf [AR]}, p.106).

\bigskip\noindent
{\bf 6. Example: p = 7, r = 2}. We know from Theorem 3.11 that there exist exactly two
non-isomorphic stable rank 2 bundles with trivial determinant. Let us prove it without
using topology. We use the following well-known description of the moduli space of
semi-stable rank 2 bundles over a compact Riemann surface of genus 3 (see {\bf [NR]}):

\medskip
\plainproclaim  Theorem 6.1. Let $SU_X(2)$ be the moduli space of semi-stable rank 2 bundles
with trivial determinant over a compact Riemann surface of genus $3$. Then
there is an embedding $\Phi:SU_X(2)\hookrightarrow \bbP(H^0(J^{2}(X),\calO(2\Theta))\cong
\bbP^7$, where 
$J^{2}(X)$ is the Picard variety of divisor classes of degree $2$ and $\Theta$ is the
hypersurface of effective divisor classes. For every $E\in SU_X(2)$ its image is a
divisor in $\bbP(H^0(J^{2}(X),\calO(2\Theta))) = |2\Theta|$ whose support is equal to the
set of $L\in J^2(X)$ such that $H^0(X,E\otimes L) \ne 0$.

\medskip
\plainproclaim  Lemma 6.2. Let $V_6$ be the unique irreducible $6$-dimensional representation
of
$PSL(2,\bbF_7)$. Then there are isomorphisms of $SL(2,\bbF_7)$-representations
$$V_6\cong V_6^* \cong S^2(V_-)\cong \Lambda^2(V_+),$$
$$S^2(V_6) \cong V_6\oplus S^4(V_-^*).$$

{\sl Proof.} The first isomorphism was observed already in Theorem 5.1. To see the second one use
that $S^4(V_-^*)$ is naturally included in $S^2(V_6) = S^2(S^2(V_-))$. Also $V_6$ is 
included in $S^4(V_-^*)$ as the space of quadrics containing the Veronese surface
$\nu_2(\bbP(V_-)$. It remains to compare the dimensions.

\medskip
\plainproclaim  Corollary 6.3. Let $X = X(7)$ be the Klein quartic. Then the group $G =
PSL(2,\bbF_7)$ acts naturally on $SU_X(2)$ and has exactly three fixed points represented
by the trivial bundle and two stable bundles.

{\sl Proof.} By construction the map $\Phi$ from theorem 6.1 is
$PSL(2,\bbF_7)$-invariant. So the group $SL(2,\bbF_7)$ acts linearly in
$H^0(J^{2}(X),\calO(2\Theta))$. Consider the embedding of
$X$ in $\bbP^5$ given by the linear system $|2K_X|$. Since $H^0(X,\calO(2K_X)) =
S^2(V_-^*)$ we see from Lemma 6.2 that $X$ is embedded equivariantly in $\bbP(V_6)$ where
$V_6$ is the unique irreducible 6-dimensional representation of $PSL(2,\bbF_7)$.  We
have the restriction map
$$r:H^0(J^{2}(X),\calO(2\Theta))\to H^0(\Theta,\calO_\Theta(2\Theta))$$ whose kernel is
one-dimensional and is spanned by a section with the divisor of zeroes equal to
$2\Theta$. This gives the decomposition of representations
$$H^0(J^{2}(X),\calO(2\Theta)) = \bbC\oplus H^0(\Theta,\calO_\Theta(2\Theta)).\eqno
(6.1)$$
 Now one can show that there is a canonical isomorphism of representations
$$H^0(\Theta,\calO_\Theta(2\Theta)) \cong H^0(\bbP(V_6), \calI_X(2)),$$
where $\calI_X$ is the ideal sheaf of $X$ embedded in $\bbP(V_6)$ (see {\bf [BV]}, 4.12). 
The restriction map $H^0(\bbP(V_6),\calO(2))\to H^0(X,\calO_X(4K_X))$ is
surjective (because $X\subset \bbP^5$ is projectively normal) and its kernel is
isomorphic to $H^0(\bbP(V_6), {\calI}_X(2))$. This gives an isomorphism of
representations
$$S^2(V_6)\cong H^0(\bbP(V_6), {\calI}_X(2))\oplus H^0(X,\calO_X(4K_X)).$$
Since $S^4(V_-^*) \cong H^0(X,\calO_X(4K_X))\oplus \bbC$, we obtain from Lemma 6.2
$$H^0(\bbP(V_6), \calI_X(2))\cong V_6\oplus \bbC.$$
Collecting everything together we get an isomorphism of $SL(2,\bbF_7)$-representations
$$H^0(J^{2}(X),\calO(2\Theta)) \cong \bbC\oplus
\bbC\oplus V_6.\eqno (6.2)$$
 This shows that the set of fixed points of $PGL(2,\bbF_7)$ in
$\bbP(H^0(J^{2}(X),\calO(2\Theta)))$ is equal to the line $\ell = \bbP(\bbC\oplus\bbC)$.
It remains to see that it intersects $SU_X(2)$ at 3 points. One point corresponds
to the trivial bundle and the other two to stable bundles. Let $\bbC$ be one of the
trivial one-dimensional summands in $H^0(J^{2}(X),\calO(2\Theta))$ which corresponds to
$Ker(r)$. It follows from the construction of $\Phi$ that the corresponding point in 
$|2\Theta|$ is equal to the divisor $2\Theta$ which is the value of $\Phi$ at the trivial
bundle. The map $r$ defines a projection map from $|2\Theta|\setminus\{2\Theta\}$ to the
space of  quadrics $H^0(\bbP(V_6), \calI_X(2))$. The line $\ell$ is the closure of
the fibre of this projection over the $PSL(2,\bbF_p)$-invariant quadric containing the
curve $X$. This quadric can be identified with the Grassmanian $G(2,4)$ in the Pl\"ucker
space $\bbP^5$. By Proposition 1.19 and Theorem 3.3 of {\tenbf [BV]} the intersection
$(\ell\setminus \{2\Theta\})\cap SU_X(2)$ consists of two stable bundles. One is the the
restriction to $X$ of the universal quotient bundle twisted by $\lambda^{-1}$ and another
is the dual of the restriction of the universal subbundle twisted by $\lambda$. This
proves the assertion.

\medskip
So we know how to construct the two stable $G$-invariant bundles on $X(7)$ with
trivial determinant. We embed
$X(7)$ in 
$\bbP(S^2(V_-)^*) =\bbP^5$ by $|2K_X|$. Then identify the representations 
$S^2(V_-)$ and $\Lambda^2(V_+)$, consider the Grassmanian $G(2,V_+)$ embedded by the
Pl\"ucker map, and then restrict to $X(7)$ the universal bundle and the universal
subbundle and twist them to get the trivial determinant. The ARK bundle corresponds
to the universal subbundle.

\medskip
The next lemma must be a special result of computations from {\tenbf [AR]}, pp.
101-105, however some typographical errors make it an unsuitable reference. We refer
to {\tenbf [A3]} for  the corrections and more general results. 

\medskip
\plainproclaim  Lemma 6.4. There is an isomorphism of representations of $SL(2,\bbF_7)$:
$$H^0(X(7),\lambda^k) \cong \cases V_8'&\text{ if $k = 5$}\\
V_+\oplus V_8'&\text{if $k = 7$}\\
V_+\oplus V_6'\oplus V_6'{}^*&\text{if $k = 9$}\endcases$$
where $V_8'$ and $V_6'$ are the $8$-dimensional and the $6$-dimensional irreducible
representations of the group $SL(2,\bbF_7)$ on which $-1$ does not act identically.

{\sl Proof.} By a theorem of H. Hopf (cf. {\tenbf [ACGH]}, p.108), given any linear map
$$f:A\otimes B\to C$$
where $A,B,C$ are complex linear spaces and $f$ is injective on each factor separately,
then 
$$\dim f(A\otimes B) \ge \dim A+\dim B -1.$$
We apply it to the map
$$H^0(X,\lambda^2)\otimes H^0(X,\lambda^3) = V_-^*\otimes V_+^*\to H^0(X,\lambda^5).$$
Using Table 3, Appendix 2, we find that $V_-^*\otimes V_+^* \cong V_+\oplus V_8'$.
Thus 
$H^0(X,\lambda^5)$ must contain $V_8'$ as its direct summand. By Riemann-Roch, it must
be the whole space.
Similarly, considering the map
$$H^0(X,\lambda^2)\otimes H^0(X,\lambda^5) = V_-^*\otimes V_8'\to H^0(X,\lambda^7)$$
we get that its image is of dimension $\ge 10$. Since, by Riemann-Roch, 
$H^0(X,\lambda^7)$ is of dimension $12$, and $H^0(X,\lambda^5)$ does not have vectors
invariant with respect to $SL(2,\bbF_7)$ (use that $-1$ acts non-trivially) , we obtain
that the multiplication map is surjective.  Using Table 4, Appendix 2, we find that 
$$V_-^*\otimes V_8' = V_+\oplus
V_6'\oplus V_6'{}^*\oplus V_8'.$$
This gives us two possibilities: $H^0(X,\lambda^7) = V_+\oplus V_8'$ or
$H^0(X,\lambda^7) = V_6'\oplus V_6'{}^*$. If the second case occurs, we consider
the map $H^0(X,\lambda^5)\otimes H^0(X,\lambda^7) \to H^0(X,\lambda^{12}$). Again, we
find that the image is of dimension $\ge 19$. We know the direct sum decomposition for 
the $SL(2,\bbF_7)$-module $H^0(X,\lambda^{12})$ (see Appendix 1).
It tells us that it contains exactly one one-dimensional summand. However,
Tables 3 and 4 from Appendix 2 tells us that $H^0(X,\lambda^5)\otimes H^0(X,\lambda^7) =
V_8'\otimes (V_6'\oplus V_6'{}^*)$ does not contain one-dimensional summands. Hence
the one-dimensional summand must be in the cokernel of the multiplication map. But
since the cokernel is of dimension $\le 3$ and $SL(2,\bbF_7)$ does not have
two-dimensional irreducible representations, we get a contradiction.

Finally $|\lambda^9|$ is cut out by the linear system of cubics in $\bbP(V_+)$. We know
from Table 2, Appendix 1, that $S^3(V_+)^* = V_+\oplus V_+^*\oplus V_6'\oplus
V_6'{}^*$. 
 The linear system of cubics spanned by
the polars of the unique invariant quartic in $\bbP(V_+)$  realizes the
summand $V_+$. The
$A$-curve $X(7)$ is contained in a linear system of cubics isomorphic to $V_+^*$.
This implies that $H^0(X,\lambda^9) \cong V_+\oplus V_6'\oplus V_6'{}^*$.

\medskip
\plainproclaim  Theorem 6.5. Let $E$ be the ARK bundle over $X(7)$. Then $E = N^*\otimes
\lambda^{5}$, where $N$ is the  normal bundle of  the
$A$-curve
$X(7)$ in
$\bbP(V_+) = \bbP^3$.

{\sl Proof.} By definition of the normal bundle we have the following exact 
sequences of $SL(2,\bbF_7)$-linearized bundles over $X$
$$0\to \lambda^{-2} \to T \to N \to 0.$$
$$0\to\calO_X \to V_+\otimes \lambda^3 \to T\to 0.$$
Combining them together we get an exact sequence
$$0\to F\to V_+\otimes \lambda^3 \to N \to 0,\eqno (6.3)$$
where
$$0\to \calO_X \to F\to \lambda^{-2}\to 0.\eqno (6.4)$$
First we see that 
$$\det (N) = \lambda^{12}\otimes\det F^{-1} = \lambda^{12}\otimes
\lambda^2 = \lambda^{14}.$$
Twisting (6.3) by $\lambda^{-5}$, we get
$$0\to F\otimes \lambda^{-5}\to V_+\otimes \lambda^{-2} \to N' \to 0,\eqno (6.5)$$
where $N' = N\otimes \lambda^{-5}$ with $\det N' = \lambda^4 = K_X^2$.
I claim that $N' = E^*$. Taking the exact sequence of cohomology for (6.3) and (6.5) we
get 
$$H^0(X,N') = Ker(H^1(X,F\otimes \lambda^{-5})\to H^1(X,V_+\otimes \lambda^{-2}),$$
$$0\to H^1(X,\lambda^{-5})\to H^1(X,F\otimes \lambda^{-5}) \to H^1(X,\lambda^{-7})\to
0.$$
Everything here is in the category of $SL(2,\bbF_7)$-modules. By Serre's duality
$$H^1(X,\lambda^{-5})\cong H^0(X,\lambda^7)^*, \quad H^1(X,\lambda^{-7})\cong
H^0(X,\lambda^9)^*.$$
Applying Lemma 6.4, we get 
$$H^1(X,F\otimes \lambda^{-5})\cong V_+^*\oplus V_+^*\oplus V_8'\oplus V_6'\oplus
V_6'{}^*.$$
Now, using Tables 3 and 4 from Appendix 2, we obtain
$$H^1(X,V_+\otimes \lambda^{-2}) = V_+\otimes H^1(X,\lambda^{-2}) = $$
$$V_+\otimes
H^0(X,\lambda^4)^* = V_+\otimes V_6 =  V_+^*\oplus V_8'\oplus V_6'\oplus
V_6'{}^*.$$
This gives us a summand $V_+^*$ in $H^0(X,N')$ and defines a map of vector bundles
$$\psi:V_+^*\otimes \calO_X\to N'.$$
Since the map $\Lambda^2 V_+ \to \Lambda^2 H^0(X,N')$ is nonzero and $\Lambda^2 V_+$ is irreducible, the
image of
$\psi$ is a vector subbundle of $N'$ of rank 2. The cokernel of $\psi$ is
concentrated over a finite set of points in $X(7)$ contained in the set $S$ of zeroes
of $s\wedge s'$ for some sections of $\Lambda^2(N')$. However, $\Lambda^2(N') =
K_X^2$ and hence $S$ consists of at most 8 points. Since $S$ is obviously
$G$-invariant, this is impossible. So $\psi$ is surjective. This implies that $N'$
defines an equivariant embedding of $X(7)$ in $G(2,V_+)$. The composition of this
embedding with the Pl\"ucker map is given by $|2K_X|$ and, hence is defined uniquely
(since there is only one equivariant isomorphism $S^2(V_-)^*\cong \Lambda^2(V_+)^*$).
This shows that the restriction of the universal subbundle on the Grassmannian to the
image of $X(7)$ coincides with the ARK bundle. Hence $N'{}^*$ is this bundle.

\medskip
 Since
$\lambda$  is of degree 2 and hence non-effective, and 
$H^0(X(7),\lambda^2) =  V_-$
 the condition (*) in
Theorem  5.2 is obviously satisfied. Thus the ARK-bundle is stable.

\medskip 
\plainproclaim  Theorem 6.6. Let $E$ be the ARK bundle over $X(7)$. Then there is an exact
sequences of $SL(2,\bbF_7)$-linearized bundles:
$$0\to E\otimes \lambda^{-3}\to V_-\otimes \calO_X\to \lambda^{10}\to 0,\eqno (6.6)$$
$$0\to \lambda^{-11}\to E\otimes \lambda\to \lambda^{11}\to 0.$$

{\sl Proof.} Let $N$ be the normal bundle of the $A$-curve $X(7)$. Using the proof of the
previous theorem we find the exact sequence
$$0\to F\otimes \lambda^{-2}\to V_+\otimes \lambda\to N\otimes \lambda^{-2}\to 0.\eqno
(6.7) $$
 Taking cohomology we easily get
$$H^0(X,N\otimes\lambda^{-2}) = H^1(X,F\otimes \lambda^{-2}) =
H^1(X,\lambda^{-2})\oplus H^1(X,\lambda^{-4}) =$$
$$ H^0(X,\lambda^{4})^*\oplus H^0(X,\lambda^{6})^* = V_6\oplus V_3^*\oplus V_7.$$
This defines a map $V_+^*\to N\otimes\lambda^{-2}$. Since $\det N\otimes\lambda^{-2} =
\lambda^{10}$, and there are no $G$-invariant subsets in $X(7)$ of cardinality 
$\le {\tenrm deg} \lambda ^{10} = 20$,
using the same argument as in the proof of Theorem 6.5, we conclude that this map is
surjective. Thus we have an exact sequence
$$0\to \lambda^{-10}\to V_-^*\otimes\calO_X\to N\otimes\lambda^{-2}\to 0.$$
After dualizing and using  Theorem 6.5, we obtain the sequence (6.6).

To get the second exact sequence, we twist (6.7) by $\lambda^6$ to obtain
$$0\to F\otimes \lambda^{4}\to V_+\otimes \lambda^7\to N\otimes \lambda^{4}\to 0.\eqno
(6.8) $$
The exact sequence
$$0\to \lambda^4\to F\otimes \lambda^{4}\to \lambda^{2}\to 0$$
gives 
$$H^1(X,F\otimes \lambda^{4}) = H^1(X,\lambda^2) =\bbC,\quad H^0(X,F\otimes
\lambda^{4}) = V_6\oplus V_-^*.$$
Since $H^1(X,V_+\otimes \lambda^7) = 0$, we obtain that 
$H^0(X,N\otimes \lambda^{4})$ is mapped surjectively onto $H^1(X,F\otimes \lambda^{4})$,
and hence contains a $G$-invariant section. This section defines a non-trivial (and
hence injective) map 
$\calO_X\to N\otimes \lambda^{4}$ and  hence an injective  map 
$$
\lambda^{-11} \to N\otimes \lambda^{-7} = E\otimes \lambda.$$
The cokernel of this map does not have torsion since otherwise its support will be a
$G$-invariant subset of cardinality $\le {\tenrm deg} \lambda^{11} = 22$. The smallest
cardinality of a $G$-invariant set on $X(7)$ is 24. Thus the cokernel is isomorphic to
$\lambda^{11}$.

\medskip\noindent
{\bf Remark 6.7} The
exact sequence (6.6) can be defined as follows. We have a linear system of curves of
degree 5 which are polars of a unique
$G$-invariant sextic in
$\bbP(V_-)$. This defines a map $V_-\to H^0(X,\calO_X(5K_X)  = H^0(X,\lambda^{10})$ and
also a surjective map
$V_-\otimes \calO_X\to \lambda^{10}$.

\medskip\noindent
{\bf Remark 6.8} Let $E$ be the ARK bundle over $X(7)$. The ruled surface $\bbP(E)$ can
be projected to $\bbP(V_+)$ with the image equal to the tri-secant scroll $S$ of the
$A$-curve $X(7)$.  Recall that any even-theta characteristic
$\theta$ on a canonical curve $X$ of genus 3 defines a (3-3)-correspondence $R$ on $X$
(see, for example, {\tenbf [DK]}):
$$R = \{(x,y)\in X\times X:|\theta+x-y|\ne \emptyset\}.$$
The set $R(x)$ consists of points in the unique positive divisor equivalent to
$\theta+x$. Let $X'$ be the image of $X$ in $\bbP^3$ under the map given by the
complete linear system $|3\theta|$. Since 
$|3\theta -(\theta+x)| = |2\theta-x| = |K_X-x|$
is a pencil, we obtain that the image of $R(x)$ lies on a line $<R(x)>$. This is a
tri-secant line of $X'$. When $x$ runs over the set of points of $X$, the tri-secant
lines form a scroll $S$ containing $X'$ as a singular curve of multiplicity 3. The map
$x\to <R(x)>$ defines an embedding of $X$ in $G(2,4)$. This
scroll also can be described as follows (see {\tenbf [SR]}). The linear system of
cubics through
$X'$ is of dimension $3$ and defines a birational transformation of $\bbP^3$ to 
$\bbP^3$. It factors through the blow-up $Y\to \bbP^3$ of $X'$ and blow-down of the
proper inverse transform $S'$ of $S$ to $X'$ in another copy of $\bbP^3$. This shows
that $S'$ is isomorphic to the projectivization of the
normal bundle of $X'$. In our case the tri-secant scroll $S$ is a surface in 
$\bbP(V_+)$ defined by an invariant polynomial of degree 8 (see {\tenbf [E1]}). 
Also in our case the correspondence $R$ on $X(7)$ is a modular correspondence $T_2$. This was
discovered by F. Klein {\tenbf [KF]} (see also {\tenbf [A2]}).

\medskip
What is the second stable $G$-invariant bundle with trivial determinant over $X(7)$? It is 
very well-known.
It can be described, for example, as follows.  Embed $X(7)$ into the Jacobian $Jac^1(X)$ and take the
normal bundle. In other words, consider a natural bijective map $V_-^*\to H^0(X,\calO_X(K_X))$. It defines
an exact sequence
$$0\to E'\to V_-^*\otimes \calO_X\to \lambda^2\to 0.\eqno (6.9)$$
Then $ E = E'\otimes\lambda$ is a $G$-invariant rank 2 bundle with trivial
determinant. 

Another way to see it is to restrict the tangent bundle of $\bbP(V_-)$ to the Klein quartic. This bundle is
$E(\lambda^3)$.

\medskip
\plainproclaim  Theorem 6.9. Let $E = E'\otimes\lambda$, where $E'$ is defined by the exact
sequence (6.9). Then $E$ is a stable rank 2 $G$-invariant bundle over $X(7)$ which
admits a non-split $G$-invariant extension
$$0\to \lambda ^{-5}\to E\to \lambda^5 \to 0.\eqno (6.10)$$

{\sl Proof.} Twisting (6.6) by $\lambda^6$ we obtain the exact sequence
$$0\to E\otimes\lambda^5\to V_-^*\otimes \lambda^6\to \lambda^8\to 0.\eqno (6.11)$$
This gives
$$H^0(X,E\otimes\lambda^5) = Ker(V_-^*\otimes H^0(X,\lambda^6)\to H^0(X,\lambda^8))=$$
$$Ker(V_-^*\otimes (V_-\oplus V_7)\to V_6\oplus V_8),$$
where $V_7$ and $V_8$ are the 6-dimensional and the 8-dimensional irreducible
representations of the group $PSL(2,\bbF_7)$. Since $V_-^*\otimes V_-$ contains a
one-dimensional summand, we find a
$G$-invariant section of $E\otimes\lambda^5$. Arguing as in the proof of Theorem 6.5, we
get the extension (6.10). 

Now let us check the stability of $E$. Assume that $E$ is not stable. Then it must
contain a $G$-invariant destabilizing subbundle $\lambda^a$ for some $a\ge 0$. Twisting
(6.9) by
$\lambda^{1-a}$  we obtain an extension
$$0\to E\otimes\lambda^{-a} \to V_-\otimes \lambda^{1-a}\to \lambda^{3-a}\to 0,$$
where $H^0(X,E\otimes\lambda^{-a}) \ne 0$. However $H^0(X,V_-\otimes \lambda^{1-a})$
does not have
$SL(2,\bbF_7)$-invariant sections when $a \ge 0$. This proves that $E$ is stable.

\medskip
We denote by $E_{(-5,5)}$ and $E_{(-11,11)}$ the two non-isomorphic stable $G$-invariant
rank 2 bundles with trivial determinant over $X(7)$ with corresponding exponential
sequences. Notice that the exponential sequences agree with Theorem 4.2 and Corollary
4.3.
 As we know from the proof of Theorem 6.1 the bundles $E_{(-5,5)}$ and $E_{(-11,11)}$ are
dual in the following sense: there exists an exact sequence
$$0\to E_{(-11,11)}\otimes\lambda^{-2}\to V_+\otimes\calO_X\to
E_{(-5,5)}\otimes\lambda^2\to 0
\eqno (6.12)$$
We have seen already that $E_{(-11,11)}\otimes\lambda^{-2}$ is the ARK bundle and its
projectivization is a non-singular model of the tri-secant scroll of the $A$-curve
$X(7)$. The projectivization of the bundle $E_{(-5,5)}\otimes\lambda^{-2} =
(E_{(-5,5)}\otimes\lambda^2)^*$ is a nonsingular model of a scroll $S^*$ in the dual
space
$V_+^*$. We can change the roles of $V_+$ and $V_+^*$ by changing the action of
$PSL(2,\bbF_7)$ on $X(7)$ via an outer automorphism of the group. Then we consider
$S^*$ as a scroll of degree  8 in the $A$-space invariant with respect to
$PSL(2,\bbF_7)$. It is equal to the
Hessian of the quartic $G$-invariant surface in $\bbP(V_+)$ (see {\bf [E1]}).  

\medskip\noindent
{\bf Remark 6.10} There are some  natural non-stable $G$-invariant rank 2 bundles over
$X(7)$. For example, the bundle $F$ which is defined by the (non-split) exact sequence
(6.4).  Tensoring with $\lambda$ we see that $F$ is unstable. A sequence of
exponents of $F\otimes \lambda$ is equal to $(1,-1)$. The same bundle can be obtained by
a construction similar to the construction of
$E_{(-5,5)}$ and
$E_{(-11,11)}$. Using the polar linear system of the Klein quartic we realize $V_3$ as a
submodule of
$S^3(V_3)^*$ and  obtain an exact sequence
$$0\to E' \to V_-\otimes\calO_X\to \lambda^6 \to 0.$$
The bundle $E = E'\otimes \lambda^{3}$ is a $G$-invariant rank 2 bundle with trivial
determinant. To see that it is unstable we use that
$$H^0(X,E\otimes\lambda^{-1}) = Ker (V_-\otimes H^0(X,\lambda^2)\to H^0(X,\lambda^8)) =$$
$$Ker(V_-\otimes V_-^*\to V_6\oplus V_8) = Ker(\bbC\oplus V_8\to V_6\oplus V_8)
\supset \bbC.$$
This shows that $E$ contains $\lambda$ as its subbundle. It is easy to see that the
quotient is the line bundle $\lambda^{-1}$ and the corresponding extension does not
split. Since ${\tenrm Ext}^1(\lambda^{-1},\lambda) = H^1(X,\lambda^2) =\bbC$, we obtain
that $E$ is isomorphic to $F\otimes\lambda$.

 \bigskip\noindent
{\bf 7. Example: p = 7, r = 3}. We use the following:

\medskip
\plainproclaim  Theorem 7.1 ({\tenbf [Bo]}). Let $p = 6n\pm 1$. The number of
irreducible 3-dimensional unitary representation of the Brieskorn sphere
$\Sigma(2,3,p)$ is equal to $3n^2\pm n$. All these representations are trivial on the
center of
$\pi_1(\Sigma(2,3,p))$.

\medskip
Applying Theorem 3.7, we obtain 

\medskip
\plainproclaim  Corollary 7.2. Let $p = 6n\pm 1$ and $G = PSL(2,\bbF_p)$. There are exactly
$3n^2\pm n$ non-isomorphic $G$-stable rank 3 bundles with trivial determinant over
$X(p)$.  Moreover each such a bundle admits a unique $G$-linearization. If $p\ne 7$ each $G$-stable bundle
is stable.

{\sl Proof.} Only the last assertion does not follow immediately from Theorem 7.1 and 
Theorem 3.7. Let us
prove it. Suppose that $E$ is $G$-stable but not stable. Let 
$0\subset E_1\subset \ldots
\subset E$ be the associated Harder-Narasimhan filtration (see {\tenbf [Se]}, p.18).
The corresponding graded bundle 
$Gr(E) =\oplus_i E_i/E_{i-1}$ is defined uniquely by $E$ and hence $G$-invariant. Since the slope of each
graded piece is equal to zero, and we know that there are no $G$-invariant line bundles of degree zero
except the trivial one, we obtain that each graded piece has trivial determinant. Assume that one of the
graded pieces is of rank 2. Then it or the other piece is a $G$-invariant subbundle of $E$. Since $E$ is
$G$-stable this is impossible. Thus all graded pieces are isomorphic to the trivial line bundle. This
implies that $Gr(E) = \calO_X^3$ is the trivial $G$-stable bundle and hence its
$SL(2,\bbF_p)$-linearization is defined by a 3-dimensional irreducible representation of $SL(2,\bbF_p)$.
It exists only for $p =7$.

\medskip
Let us assume $p = 7$. We need to exhibit four  rank 3 $G$-stable vector bundles with
trivial determinant over
$X(7)$. This is easy. First of all we take the two rank 2 bundles $E_{(-5,5)}$ and 
$E_{(-11,11)}$ and consider their second symmetric powers $S^2(E_{(-5,5)}),
S^2(E_{(-11,11)})$. The other two are obtained by considering the trivial bundles 
$V_-\otimes\calO_X, V_+\otimes\calO_X$ with linearizations defined by two irreducible
3-dimensional representations of $G$. It is easy to compute their sequences of
exponents $(a_1,a_2,a_3)$:

\medskip
\plainproclaim  Theorem 7.3. We have
$$(a_1,a_2,a_3) =\cases (-10,0,10)&\text{if $E =  S^2(E_{(-5,5)})$}\\
(-22,0,22)&\text{if $E = S^2(E_{(-11,11)})$}\\
(-2,-4,6)&\text{if $E = V_-\otimes\calO_X$}\\
(-6,4,2)&\text{if $E = V_-^*\otimes\calO_X$}\endcases.$$

{\sl Proof.} The first two sequences can be immediately computed from the known sequences
of exponents of the rank 2 bundles. To compute the third sequence we use that 
$H^0(X,V_-\otimes \lambda^2) = V_-\otimes V_-^* $ contains a trivial summand. This
implies that $E$ contains a $G$-linearized subbundle isomorphic to $\lambda^{-2}$. Using
exact sequence (6.9) we see that the quotient $F$ is isomorhic to
$E_{(-5,5)}\otimes\lambda$. Now we use that $(-4,6)$ is a sequence of exponents for
$E_{(-5,5)}\otimes\lambda$. From the duality we obtain that $F^*$ is a $G$-linearized
subbundle of $V_-^*\otimes\calO_X$ and the quotient is $\lambda^2$. This gives the
fourth sequence.

\medskip\noindent
{\bf Remark 7.4} The first two bundles correspond to unitary representations of
the group $\pi_1(\Sigma(2,3,7))$ which arise from a representation
$\rho:\pi_1\to SO(3)\subset SU(3)$. The remaining two bundles correspoond to
''additional'' ({\tenbf [Bo]}, p.211) irreducible representations. A new information
is that the additional representations factor through an irreducible representation of
$PSL(2,\bbF_7)$. This was verified directly by H. Boden.

\medskip
We have a canonical exact sequence corresponding to the embedding of $X(7)$ as the
$A$-curve:
$$0\to E'\to V_+^*\otimes\calO_X\to \lambda^3\to 0 \eqno (6.13)$$
Twisting by $\lambda$ we obtain a rank 3 $G$-invariant bundle $E$ with trivial
determinant. One can show that $E\cong S^2(E_{(-11,11)})$. 

Similarly, the exact sequence
$$0\to E'\to V_+^*\otimes\calO_X\to \lambda^9\to 0$$
defined by the polar linear system of the $SL(2,\bbF_7)$-invariant quartic surface in the
$A$-space defines the $G$-invariant rank 3 bundle $E = E'(\lambda^3)$ with trivial
determinant. Once can show that it is isomorphic to $S^2(E_{(-5,5)})$.

\bigskip\noindent
{\bf 8. Example: p = 11, r = 2}. By Theorem 3.11 we expect to find four
non-isomorphic $PSL(2,\bbF_{11})$-invariant stable bundles of rank 2 with trivial determinant. 
Here
we have of course the ARK bundle which is stable if one checks that
$V_+^*\not\subset H^0(X(11),\lambda^{a})$ for
$a\le 4$. Assume $H^0(X(11),\lambda^{4})$ contains $V_+$. Since it contains already
$V_-^*$ we would have $\dim |\lambda^4|
\ge 11$. This contradicts the Clifford theorem (see {\bf [ACGH]}). The same theorem implies that 
$V_+^*\not\subset H^0(X(11),\lambda^{a})$ for $a\le 2$. So, it remains to verify that 
$V_+^*\not\subset H^0(X(11),\lambda^{3})$. I do not know how to do it.

Another potential candidate is the vector bundle defined by using the fact that
$z$-curve
$X(11)$ parametrizes polar quadrics of corank $2$ of the invariant cubic hypersurface
$W$ (see Example 2.9). This defines a bundle with determinant $\lambda^{6}$ which
embeds $X(11)$ in $G(2,V_-)$. The corresponding ruled surface in $\bbP(V_-)$ is the four-secant scroll 
of the $z$-curve $X(11)$ of degree 30 (see {\bf [E2]}). Similar to the tri-secant scroll of $X(7)$ it is
defined by a modular $(4,4)$-correspondence on $X(11)$ (see {\bf [A2]}).

To introduce the third candidate, we use that the cubic hypersurface $W$ admits a
$G$-invariant representation as the Pfaffian hypersurface (see {\tenbf [AR]}, p. 164):
$$v^2w+w^2x+x^2y+y^2z+z^2v = {\tenrm Pf}\pmatrix 0&v&w&x&y&z\\
-v&0&0&z&-x&0\\
-w&0&0&0&v&-y\\
-x&-z&0&0&0&w\\
-y&x&-v&0&0&0\\
-z&0&y&-w&0&0\endpmatrix.\eqno (8.1)$$
This representation is obtained by considering a linear map $V_-\to S^2(V_-^*)$ defined by the polar
linear system of $W$  and then identifying the
representations $S^2(V_-^*)$ and
$\Lambda^2(V_+)^*$ (see Theorem 5.1). The cubic hypersurface $W$ is equal to the
pre-image of the cubic hypersurface in 
$\bbP(\Lambda^2(V_+)^*)$ which coincides with the chordal variety $C$ of the Grassmanian
$G(2,V_+^*)$. The latter carries a canonical rank 2 bundle whose fibre over a point
$t\in C$ is equal to the null-space $L_t\subset V_+$ of the corresponding skew-symmetric
matrix.  To get a bundle over $X(11)$ we use the decomposition of
$PSL(2,\bbF_{11})$-representations
$$H^0(X,K_X) = V_-\oplus V_{10}\oplus V_{11},\eqno (8.2)$$
where $V_{10}$ and $V_{11}$ are 10-dimensional and $11$-dimensional irreducible
representations of $PSL(2,\bbF_{11})$ (see {\tenbf [He]}). 
This decomposition allows us to project $X(11)$ to $\bbP(V_-)$ as a curve of degree $50$.
It turns out that the image of the projection is contained in the cubic $W$. This result
is due to F. Klein {\tenbf [Kl]} and is reproduced by A. Adler (see {\tenbf [AR]},
Appendix 3). Other proof is due to M. Gross and S. Popescu {\tenbf [GP]} who give an
interpretation of the cubic $W$ as a compactification of the moduli space of abelian
surfaces with polarizaton of type $(1,11)$. This allows us to restrict the bundle
$E$ to obtain a
$G$-invariant bundle over
$X(11)$. This bundles embeds $X(11)$ in $G(2,V_+)$ by the linear system of quadrics spanned by 
the Pfaffians
of order four principal submatrices of the skew-symmetric matrix from (8.1). So, the determinant of the
bundle is equal to $\lambda^8$.

Finally one may try to consider the normal bundle of $X(11)$ in the cubic $W$.	
I do not know yet whether any of these four bundles is stable, nor do I know their sequence of
exponents. 
\vfill\eject

\bigskip\noindent
DEPARTMENT OF MATHEMATICS, UNIVERSITY OF MICHIGAN, 
ANN ARBOR, MI 48109

\smallskip\noindent
e-mail:idolga\@math.lsa.umich.edu
\vfill\eject

\centerline{{\bf Appendix 1. Decompositions of $S^n(V_-)$ and $S^n(V_+)$ for $p = 7$}}

\medskip
It follows from the character table of the groups $SL(2,\bbF_p)$ (see, for example, {\bf [D]}) that
$SL(2,\bbF_7)$ has eleven  non-isomorphic irreducible representations. In the following we denote by
$V_k$ the irreducible representation of
$SL(2,\bbF_7)$ of dimension
$k$ and by
$V_k'$ another representation of the same dimension which does not factor through
$PSL(2,\bbF_{7})$. We assume that 
$$V_- =
V_3, \quad V_+ = V_4.$$  

\medskip
The following generating function was computed in {\bf [BI]}):
$$Q_V(t) =\sum_{t=0}^\infty \dim_\bbC {\tenrm Hom}^{PSL(2,\bbF_7)}(V,S^n(V_-))t^n.$$ 
We have
$$Q_{V_1}(t) = {1+t^{21}\over (1-t^4)(1-t^6)(1-t^{14})},$$
$$Q_{V_3}(t) = {1+t^7+t^{11}+t^{13}\over (1-t^4)(1-t^7)(1-t^8)},$$
$$Q_{V_3^*}(t) = {t^3+t^5+t^9+t^{15}\over (1-t^4)(1-t^7)(1-t^8)},$$
$$Q_{V_6}(t) = {t^2+t^8\over (1-t^2)(1-t^4)(1-t^7)},$$
$$Q_{V_7}(t) = {t^3\over (1-t^2)(1-t^3)(1-t^4)},$$
$$Q_{V_8}(t) = {t^4\over (1-t)(1-t^3)(1-t^7)}.$$ 

Using the similar arguments one can compute the
 generating function
$$P_V(t) =\sum_{t=0}^\infty \dim_{\bbC} {\tenrm Hom}^{SL(2,\bbF_7)}(V,S^n(V_+))t^n.$$

We have
$$P_{V_1}(t) = {1+t^8+t^{10}+t^{12}+t^{16}+t^{18}+t^{20}+t^{28}\over
(1-t^4)(1-t^6)(1-t^8)(1-t^{14})},$$
$$P_{V_3}(t) =
{t^2-t^4+t^6+2t^8+t^{12}+2t^{14}+2t^{18}\over
(1-t^2)(1-t^4)(1-t^8)(1-t^{14})},$$
$$P_{V_3^*}(t) = {2t^6+2t^{10}+t^{12}+2t^{16}+t^{18}-t^{20}+t^{22}\over
(1-t^2)(1-t^4)(1-t^8)(1-t^{14})},$$
$$P_{V_6}(t) = {2t^4+t^8+2t^{10}+t^{12}+2t^{16}\over(1-t^2)(1-t^4)^2(1-t^{14})},$$
$$P_{V_7}(t) = {t^2+t^4+2t^6+2t^{10}+t^{12}+t^{14}\over
(1-t^2)(1-t^4)(1-t^6)(1-t^8)},$$
$$P_{V_8}(t) = {t^4+t^6+2t^8+2t^{12}+t^{14}+t^{16}\over
(1-t^2)^2(1-t^6)(1-t^{14})},$$
$$P_{V_4}(t) = {t^3+t^7-t^9+2t^{11}+t^{15}-t^{17}+t^{19}\over
(1-t^2)^2(1-t^6)(1-t^{14})},$$
$$P_{V_4^*}(t) = {t-t^3+t^5+2t^9-t^{11}+t^{13}+t^{17}\over
(1-t^2)^2(1-t^6)(1-t^{14})},$$
$$P_{V_6'}(t) = P_{V_6'{}^*}(t) = {t^3+2t^9+t^{15}\over
(1-t^2)^2(1-t^4)(1-t^{14})},$$
$$P_{V_8'}(t) = {2t^5+t^7+t^9+t^{11}+t^{13}+2t^{15}\over
(1-t^2)^2(1-t^6)(1-t^{14})}.$$

\medskip
The generating functions $P_V(t)$ and $Q_V(t)$ allows one, in principle, decompose
any symmetric power $S^n(V_-^*)$ or $S^n(V_+^*)$ in irreducible representations of
$SL(2,\bbF_7)$. We give a few examples:
\medskip
\valign{&\hbox{\strut\quad#}\cr
$S^2(V_3)^*$&
$S^3(V_3)^*$&
$S^4(V_3)^*$&
$S^5(V_3)^*$&
$S^6(V_3)^*$&
$S^7(V_3)^*$&
$S^8(V_3)^*$&
$S^9(V_3)^*$&
$S^{10}(V_3)^*$&
$S^{11}(V_3)^*$\cr
= &= &= &= &=&= &= &= &= &= \cr
$V_6$&$V_3+V_7$&$V_1+V_6+V_8$&$V_3+V_3^*+V_7+V_8$&$V_1+2\cdot V_6+V_7+V_8$ &$V_3+V_3^*+2\cdot V_7+2\cdot
V_8$&$V_1+V_3^*+3\cdot V_6+V_7+2\cdot V_8$&$2\cdot V_3+2\cdot V_3^*+V_6+3\cdot V_7+2\cdot
V_8$&$V_1+V_3+4\cdot V_6+2\cdot V_7+3\cdot V_8$&$2\cdot V_3+2\cdot V_3^*+V_6+4\cdot V_7+4\cdot V_8$\cr}
\medskip

\centerline{Table 1: Decomposition of $S^n(V_3^*)$}

\medskip
\valign{&\hbox{\strut\quad#}\cr
$S^2(V_4)^*$&
$S^3(V_4)^*$&
$S^4(V_4)^*$&
$S^5(V_4)^*$&
$S^6(V_4)^*$\cr
= &= &= &= &= \cr
$V_3+V_7$&$V_4+V_{4}^*+V_{6'}+V_{6'}^*$&$V_1+2\cdot V_6+2\cdot V_7+V_8$&$2\cdot V_4+2\cdot
V_{4}^*+2\cdot V_{6'}+2\cdot V_{6'}^*+2\cdot V_{8'}$&$V_1+2\cdot V_3+2\cdot V_3^*+2\cdot V_6+5\cdot
V_7+3\cdot V_8$\cr}

\medskip

\centerline{Table 2: Decomposition of $S^n(V_4^*)$}
\medskip
To deduce from this the decompositions for $H^0(X(7),\lambda^{2n}) =
H^0(X(7),K_{X(7)}^{n})$ we have to use that
$$\sum_{n=0}^\infty\dim H^0(X(7),\lambda^{2n})t^n = (1-t^4)\sum_{n=0}^\infty\dim
S^n(V_3)^*t^n.$$

\vfill\eject
\centerline{{\bf Appendix 2. Tables for tensor products of reprsentations of  
SL(2,$\bbF_7)$}}
\smallskip
We use the notation from Appendix 1. For brevity we skip $V$ in the notation $V_n$. The
following tables give the decompositions for the tensor products of irreducible
representations of
$SL(2,{\bbF}_7)$.
\bigskip
\valign{&\hbox{\strut\quad#}\cr
{}&
$V_3$&$V_3^*$&$V_6$&$V_7$&$V_8$&$V_4$&$V_4^*$&$V_6'$&$V_6'{}^*$&$V_8'$\cr
$V_3$&$3^*+6$&$1+8$&$3^*+7+8$&$6+7+8$&$3+6+7+8$&$
4^*+8'$&$6'+6'{}^*$&$4+6'+8'$&$4+6'{}^*+8'$&$4^*+6'+6'{}^*+8'$\cr
$V_3^*$&$1+8$&$3+6$&$3+7+8$&$6+7+8$&$3^*+6+7+8$&$6'+6'{}^*$&$4+8'$&$4^*+6'+8'$&
$4^*+6'{}^*+8'$&$4+6'+6'{}^*+8'$\cr
$V_6$&$3^*+7+8$&$3+7+8$&$1+2\cdot 6+7+2\cdot 8$&$3+3^*+6+2\cdot
7+2\cdot 8$&$3+3^*+2\cdot 6+2\cdot 7+2\cdot
8$&$4^*+6'+6'{}^*+8'$&$4+6'+6'{}^*+8'$&$4+4^*+6'+6'{}^*+2\cdot
8'$&$4+4^*+6'+6'{}^*+2\cdot 8'$&$4+4^*+2\cdot 6'+2\cdot 6'{}^*+2\cdot 8'$\cr}

\bigskip
\centerline{Table 1}

\bigskip
\valign{&\hbox{\strut\quad#}\cr
{}&
$V_3$&$V_3^*$&$V_6$&$V_7$&$V_8$&$V_4$&$V_4^*$&$V_6'$&$V_6'{}^*$&$V_8'$\cr
$V_7$&$6+7+8$&$6+7+8$&$3+3^*+6+2\cdot
7+2\cdot 8$&$1+3+3^*+2\cdot 6+2\cdot 7+2\cdot 8$&$3+3^*+2\cdot
6+2\cdot 7+3\cdot 8$&$4+4^*+6'+6'{}^*+8'$&$4+4^*+6'+6'{}^*+8'$&$4+4^*+2\cdot
6'+6'{}^*+2\cdot 8'$&$4+4^*+6'+ 2\cdot
6'{}^*+2\cdot 8'$&$4+4^*+2\cdot 6'+2\cdot 6'{}^*+3\cdot 8'$\cr
$V_8$&$3+6+7+8$&$3^*+6+7+8$&$3+3^*+2\cdot 6+2\cdot 7+2\cdot
8$&$3+3^*+2\cdot 6+2\cdot 7+3\cdot 8$&$1+3+3^*+2\cdot 6+3\cdot
7+3\cdot 8$&
$4+6'+6'{}^*+2\cdot 8'$&$4^*+6'+6'{}^*+2\cdot 8'$&$4+4^*+2\cdot 6'+2\cdot 6'{}^*+2\cdot
8'$&$4+4^*+2\cdot 6'+2\cdot 6'{}^*+2\cdot 8'$&$2\cdot 4+2\cdot 4^*+2\cdot 6'+2\cdot
6'{}^*+3\cdot 8'$\cr}

\bigskip
\centerline{Table 2}

\bigskip
\valign{&\hbox{\strut\quad#}\cr
{}&
$V_3$&$V_3^*$&$V_6$&$V_7$&$V_8$&$V_4$&$V_4^*$&$V_6'$&$V_6'{}^*$&$V_8'$\cr
$V_4$&$
4^*+8'$&$6'+6'{}^*$&$4^*+6'+6'{}^*+8'$&$4+4^*+6'+6'{}^*+8'$&$4+6'+6'{}^*+2\cdot 8'$&$
3^*+6+7$&$
1+7+8$&$3+6+7+8$&$3+6+7+8$&$3^*+6+7+2\cdot 8$\cr
$V_4^*$&$6'+6'{}^*$&$
4+8'$&$4+6'+6'{}^*+8'$&$4+4^*+6'+6'{}^*+8'$&$4^*+6'+6'{}^*+2\cdot 8'$&$
1+7+8$&$
3+6+7$&$3^*+6+7+8$&$3^*+6+7+8$&$3+6+7+2\cdot 8$\cr
$V_6'$&$4+6'+8'$&$4^*+6'+8'$&$4+4^*+6'+6'{}^*+2\cdot 8'$&$4+4^*+2\cdot 6'+6'{}^*+2\cdot
8'$&$4+4^*+2\cdot 6'+2\cdot 6'{}^*+2\cdot 8'$&
$
3+6+7+8$&$3^*+6+7+8$&$6+2\cdot 7+2\cdot 8$&$1+3+3^*+6+7+2\cdot 8$&$3+3^*+2\cdot 6+2\cdot 
7+2\cdot 8$\cr}

\bigskip
\centerline{Table 3}

\vfill\eject

\bigskip
\valign{&\hbox{\strut\quad#}\cr
{}&
$V_3$&$V_3^*$&$V_6$&$V_7$&$V_8$&$V_4$&$V_4^*$&$V_6'$&$V_6'{}^*$&$V_8'$\cr
$V_6'{}^*$&$4+6'+8'$&$4^*+6'+8'$&$4+4^*+6'+6'{}^*+2\cdot 8'$&$4+4^*+6'+2\cdot
6'{}^*+2\cdot 8'$&$4+4^*+2\cdot 6'+2\cdot 6'{}^*+2\cdot 8'$&
$
3^*+6+7+8$&$3+6+7+8$&$6+2\cdot 7+2\cdot 8$&$1+3+3^*+6+7+2\cdot 8$&$3+3^*+2\cdot 6+2\cdot 
7+2\cdot 8$\cr
$V_8'*$&$4^*+6'+6'{}^*+8'$&$4+6'+6'{}^*+8'$&$4+4^*+2\cdot 6'+2\cdot
6'{}^*+2\cdot 8'$&$4+4^*+2\cdot 6'+2\cdot 6'{}^*+3\cdot 8'$&$2\cdot 4+2\cdot 4^*+2\cdot
6'+2\cdot 6'{}^*+3\cdot 8'$&$3^*+6+7+2\cdot 8$&$3+6+7+2\cdot 8$&$3+3^*+2\cdot 6+2\cdot
7+2\cdot 8$&$3+3^*+2\cdot 6+2\cdot 7+2\cdot 8$&$3+3^*+2\cdot 6+3\cdot 7+3\cdot 8$\cr }

\bigskip
\centerline{Table 4}

\bye